\documentclass{amsart} 

\newdimen\theight
\def \Column{%
             \vadjust{\setbox0=\hbox{\sevenrm\quad\quad tcol}%
             \theight=\ht0
             \advance\theight by \dp0    \advance\theight by \lineskip
             \kern -\theight \vbox to \theight{\rightline{\rlap{\box0}}%
             \vss}%
             }}%

\catcode`\@=11
\def\qed{\ifhmode\unskip\nobreak\fi\ifmmode\ifinner\else\hskip5\p@\fi\fi
 \hbox{\hskip5\p@\vrule width4\p@ height6\p@ depth1.5\p@\hskip\p@}}
\catcode`@=12 

\newcount\notenumber
\def\clearnotenumber{\notenumber=0}
\def\note{\global\advance\notenumber by 1
 \footnote{$^{\the\notenumber}$}}

\clearnotenumber

\def\Ga{{\relax\ifmmode\frak a\else$ \frak a $\fi}}
\def\ga{{\relax\ifmmode\frak a\else$ \frak a $\fi}}
\def\Gb{{\relax\ifmmode\frak b\else$ \frak b $\fi}}
\def\Gc{{\relax\ifmmode\frak c\else$ \frak c $\fi}}
\def\Gd{{\relax\ifmmode\frak d\else$ \frak d $\fi}}
\def\Ge{{\relax\ifmmode\frak e\else$ \frak e $\fi}}
\def\GA{{\relax\ifmmode\frak A\else$ \frak A $\fi}}
\def\GC{{\relax\ifmmode\frak C\else$ \frak C $\fi}}
\def\gC{{\relax\ifmmode\frak C\else$ \frak C $\fi}}
\def\GI{{\relax\ifmmode\frak I\else$ \frak I $\fi}}

\font\msym=msbm10
\font\msyms=msbm7
\def\smallbox#1{\leavevmode\thinspace\hbox{\vrule\vtop{\vbox
   {\hrule\kern1pt\hbox{\vphantom{\tt/}\thinspace{\tt#1}\thinspace}}
   \kern1pt\hrule}\vrule}\marginpar{!!}\thinspace}

\def\note{\footnote} 

\makeatletter
\def\PROOF{\par\noindent{\it Proof:} \@ifnextchar:{\@gobbleignore}{}}
\def\@gobbleignore#1{\ignorespaces}
\makeatother

   \def\nofork^#1_#2{\putforkinmargin
        \unionstick^{\textstyle #1}_{\textstyle #2}}

\newbox\noforkbox \newdimen\forklinewidth
\setbox0\hbox{$\textstyle\bigcup$}
\setbox1\hbox to \wd0{\hfil\vrule width \forklinewidth depth \dp0
                        height \ht0 \hfil}
\setbox\noforkbox\hbox{\box1\box0\relax}
\def\unionstick{\mathop{\copy\noforkbox}\limits}
\def\nonfork#1#2_#3{#1\unionstick_{\textstyle #3}#2}
\def\nonforkin#1#2_#3^#4{#1\unionstick_{\textstyle #3}^{\textstyle #4}#2}

\def\cal{\mathcal}

\def\opname#1 {\expandafter\def\csname #1\endcsname
		{\ifmmode\mathop{\mathrm {#1}}\else {\rm #1}\fi}}

\def\xname#1 {\expandafter
  \def\csname#1\endcsname{{\ifmmode\mathord{\mathrm {#1}}\else {\rm #1}\fi}}}

\def\textspace{\ifmmode\else\ \fi}
\makeatletter
\def\mathname#1 {\expandafter
  \def\csname#1\endcsname{{\ifmmode\mathord{\mathrm{#1}}\else
		\@latexerr{\commandname#1 only allowed in math mode}%
		{See comments in moredefs.tex for more help}\fi}}}

\makeatother
\def\commandname#1{\expandafter\string\csname#1\endcsname}

\def\lg{\mathop{\rm \ell g}}
\mathname Res 
\opname Dom 
\opname Max 
\opname Ker
\opname Cent 
\opname Rang 
\opname acc 
\opname accc 
\opname bd 
\opname cf
\opname cl 
\opname cov 
\opname card
\opname id 
\opname inf 
\opname lub
\opname max 
\opname min 
\opname nacc
\opname otp 
\opname pcf

\opname pp 
\opname rang 
\opname rin
\opname sup 
\opname tcf 
\opname tlim 
\xname Ord 
\xname Pr
\newcommand{\Gm}{{\Game}}
\xname Reg
\xname stp  
\xname Suc
\xname Sur
\xname ba
\xname cd
\xname is
\xname iso
\xname eub
\xname pos
\xname pr
\xname ptr
\xname tr
\xname tp
\xname trr
\xname var

\def\gothica0{{\relax\ifmmode \Ga\else $\Ga$\fi}}
\def\gothicb0{{\relax\ifmmode \Gb\else $\Gb$\fi}}
\def\gothicc0{{\relax\ifmmode \Gc\else $\Gc$\fi}}
\def\gothicd0{{\relax\ifmmode \Gd\else $\Gd$\fi}}
\def\gothice0{{\relax\ifmmode \Ge\else $\Ge$\fi}}

\newcommand{\lesdot}{\mathrel{\mathord{<}\!\!\raise 0.8
pt\hbox{$\scriptstyle\circ$}}}
\newcommand{\leseqdot}{\mathrel{\mathord{\leq}\!\!\raise 
0.8 pt\hbox{$\scriptstyle\circ$}}}

\newcommand{\conc}{{}^\frown\!}

\def\frak{\fam\frakturfam\teneuf}

 \def\nfork_#1^#2_#3{
 \mathrel{\rlap{$\displaystyle\unionstick^{\textstyle#2}_{\textstyle#3}$}
 \mathop{\phantom{\copy\noforkbox}}\nolimits_{#1}}
 }

 \newcount\skewfactor
 \def\mathunderaccent#1#2 {\let\theaccent#1\skewfactor#2
 \mathpalette\putaccentunder}
 \def\putaccentunder#1#2{\oalign{$#1#2$\crcr\hidewidth
 \vbox to.2ex{\hbox{$#1\skew\skewfactor\theaccent{}$}\vss}\hidewidth}}

 \usepackage{amssymb}

\let\mathscr\cal

\let\frak\mathfrak

\makeatletter

\newtheorem{definition}{Definition}[section]
\def\setuptheorem#1{\newtheorem{#1}[definition]{#1}}

\setuptheorem{Assumption}
\setuptheorem{Claim}
\setuptheorem{Concluding Remarks}
\setuptheorem{Conclusion}
\setuptheorem{Construction}
\setuptheorem{Context and Fact}
\setuptheorem{Context}
\setuptheorem{Convention}
\setuptheorem{Crucial Fact}
\setuptheorem{Definition and Claim}
\setuptheorem{Definition and Context}
\setuptheorem{Definition of the Game}
\setuptheorem{Definition}
\setuptheorem{Discussion}
\setuptheorem{Examples}
\setuptheorem{Example}
\setuptheorem{Exercise}
\setuptheorem{Explanation of our Intended Plan}
\setuptheorem{Explanation}
\setuptheorem{Explaning the first construction}
\setuptheorem{Fact}
\setuptheorem{Formal Description of the Translation}
\setuptheorem{Further Discussion}
\setuptheorem{Hypothesis}
\setuptheorem{Lemma}
\setuptheorem{Main Fact}
\setuptheorem{Main Lemma}
\setuptheorem{Main Definition}
\setuptheorem{Main Theorem}
\setuptheorem{Main Claim}
\setuptheorem{More Discussion}
\setuptheorem{Notation}
\setuptheorem{Observation}
\setuptheorem{Preliminaries}
\setuptheorem{Problem}
\setuptheorem{Question}
\setuptheorem{Remarks}
\setuptheorem{Remark}
\setuptheorem{Remark/Definition}
\setuptheorem{Set Theoretic Context}
\setuptheorem{Subclaim}
\setuptheorem{Subexample}
\setuptheorem{Subfact}
\setuptheorem{The Existence Theorem}
\setuptheorem{The Model Theoretic Context}
\setuptheorem{Theorem}
\setuptheorem{Thesis}

\newtheorem{the context continued}[definition]{The context continued}

\newtheorem{description of the construction}[definition]{Description of
the Construction}

\newtheorem{continuation of the proof of 5.2}[definition]{Continuation of
the proof of ???5.2\ref{what is 5.2?}}

\newtheorem{claim}[definition]{Claim}

\newtheorem{convention}[definition]{Convention}

\newtheorem{main lemma}[definition]{Main Lemma}

\def\grayon{\special{ps::  0.9 setgray }}
\def\grayoff{\special{ps:: 0 setgray }}

\setuptheorem{Comment}

\newcommand{\st}{{such that}}
\newcommand{\sq}{{sequence}}

\newcommand{\rest}{{\restriction}}

\newcommand{\power}{{cardinality}}

\newcommand{\wolog}{{without loss of generality}}
\newcommand{\Wlog}{{Without loss of generality}}
  
\newcommand{\then}{{\underline{then}}}
\newcommand{\Then}{{\underline{Then}}}
\newcommand{\there}{{\underline{there}}}

\newcommand{\when}{{\underline{when}}}

\newcommand{\uif}{{\underline{if}}}

\newcommand{\uiff}{{\underline{if and only if}}}

\setuptheorem{Main Sub Claim}

\setuptheorem{Main SubClaim}

\newcommand{\bool}{{\mathbf B}}

\newcommand{\blackw}{{\mathbf W}}

\newcommand{\idealI}{{\dot{\cal I}}}

\newcommand{\blackF}{{\mathcal F}}

\newcommand{\blackz}{{\dot \zeta}}
 
\newcommand{\bools}{{\dot s}}

\newcommand{\blackW}{{\mathbf W}}

\newcommand{\ebool}{{\dot y}}

\newcommand{\Ccc}{{c.c.c}}

\newcommand{\idealJ}{{\bf J}}

\newcommand{\rDom}{{\rm Dom}}

\newcommand{\ExKer}{{\rm ExKer}}

\newcommand{\Ext}{{\rm Ext}}

\newcommand{\rpr}{{\rm pr}}

\newcommand{\rRang}{{\rm Rang}}

\newcommand{\supp}{{\rm supp}}

\begin{document}




\def\referencemark#1#2{\ifcase #1
\or
\margin{intnl #2}%
\or
\margin{complete #2}%
\or 
\margin{incompl #2}
\fi}

\def \margin#1{%
\ifvmode\setbox0\hbox to \hsize{\hfill\rlap{\small\quad#1}}%
\ht0 0cm
\dp0 0cm
\box0\vskip-\baselineskip
\else 
             \vadjust{\setbox0=\hbox{\small\quad #1}%
             \theight=\ht0
             \advance\theight by \dp0    \advance\theight by \lineskip
             \kern -\theight \vbox to \theight{\rightline{\rlap{\box0}}%
             \vss}%
             }\fi}
\def \lmargin#1{%
\ifvmode\setbox0\hbox to \hsize{\llap{\small#1\quad}\hfill}%
\ht0 0cm
\dp0 0cm
\box0
\else 
             \vadjust{\setbox0=\hbox{\small #1\quad}%
             \theight=\ht0
             \advance\theight by \dp0    \advance\theight by \lineskip
             \kern -\theight \vbox to \theight{\leftline{\llap{\box0}}%
             \vss}%
             }\fi}

\def\aBLACKsquare{\vrule width 0.3cm height 0.3cm depth 0.2cm}

%
\def\mcol#1{\ifcase #1
%
\or
%
%
\or
%
%
\or
%
\margin{beth?}%
\or
%
%
\or
%
\margin{gothic?}
\or
%
\aBLACKsquare\margin{sup etc.}
\or
%
\margin{..bigcup..}
\or
%
\aBLACKsquare
\or
%
\margin{bigcup..}
\or
%
\margin{fixed}
\else \what\fi
}


\newcount\loopcount
\makeatletter
\def\leavefile{%
\loopcount0
\loop\catcode\loopcount14 
\advance\loopcount1 
\ifnum\loopcount<128\repeat
}
\makeatother

\everymath{\graystuff}
\def\graystuff#1${{\setbox2=\hbox{\everymath{}\relax$#1$}%
\setbox1\hbox{\vrule width\wd2 height\ht2 depth\dp2 }
\wd1=0cm
\grayon
\box1
\grayoff
\box2}$}

\def\grayon{\special{ps::  0.9 setgray }}
\def\grayoff{\special{ps:: 0 setgray }}


\everymath{}


 \everydisplay{\everymath{}}

%

\def\putforkinmargin{\margin{fork}}

\def\marginbegin#1{\margin{#1\dots}}
\def\marginend#1{\margin{\dots#1}}

\def\mcol#1{}

\def\lmargin#1{}

%


 
\makeatletter
 \def\@citex[#1]#2{%
   \let\@citea\@empty
   \@cite{\@for\@citeb:=#2\do
     {\@citea\def\@citea{,\penalty\@m\ }%
      \edef\@citeb{\expandafter\@firstofone\@citeb}%
      \if@filesw\immediate\write\@auxout{\string\citation{\@citeb}}\fi
      \@ifundefined{b@\@citeb}{\mbox{\reset@font\bfseries ?}%
        \G@refundefinedtrue
	\@latex@warning
          {Citation `\@citeb' on page \thepage \space undefined}}%
        {\hbox{\csname b@\@citeb\endcsname}}}}{#1}}

\def\citewarn[#1]#2{\@citex[#1]{#2}}
\makeatother

\newcommand{\saharon}[1]{\vv #1\marginpar{!!} \vv\ }
\newcommand{\vv}{\vrule width 5mm height 9pt depth 1pt}



\title[Existence of Endo-Rigid Boolean Algebras]{Existence of Endo-Rigid Boolean Algebras}
\author{Saharon Shelah}
\begin{abstract}
How many endomorphisms does a Boolean algebra have? Can we find Boolean
algebras with as few
endomorphisms as possible?
 Of course from any ultrafilter of the Boolean algebra
we can define an endomorphism,
and we can combine finitely many such endomorphisms
in some reasonable ways. We prove that in any cardinality
$ \lambda = \lambda ^ {\aleph_0} $ there is a Boolean algebra
with no other endomorphisms. For this we use the so  called ``black boxes'',
but in a self contained way.
We comment on how necessary  the restriction on the cardinal is.

\end{abstract}
\thanks{Publication E58. It was supposed to be Chapter I of the book
       ``Non-structure'' and probably will be if it materializes.}
\maketitle
\setcounter{section}{-1}

\bigskip

\section{Introduction}

In this paper we prove the existence of a Boolean algebra of any
cardinality $\lambda=\lambda^{\aleph_0}$ which has as few endomorphisms as
possible, in some natural sense. Note that every ultrafilter $D$ of a Boolean algebra $\bool$ induces
an endomorphism $h_D$ of $\bool$: $h_D(x)$ is $1_\bool$ for $x\in D$ and
$0_\bool$ otherwise. Also we can combine endomorphisms:
if $h_\ell$ is a homomorphism
from $\bool\restriction a_\ell$ into $\bool\restriction b_\ell$
for $\ell=1,2$ and
$a_1\cup a_2=1_\bool=b_1\cup b_2$, $a_1\cap a_2=0_\bool=b_1\cap b_2$,
{\em then\/}
there is a unique endomorphism $h$ of $\bool$ extending both $h_1$ and $h_2$, and for
any endomorphism $h$ of $\bool$ and $a\in \bool$, $h\restriction
(\bool\restriction a)$
is a homomorphism from $\bool\restriction a$ into the Boolean algebra
$\bool\restriction h(a)$.

Also if $\idealI_1,\idealI_2$ are ideals of $\bool$
satisfying $\idealI_1\cap
\idealI_2=\{0_ \bool 
\}$ and $\{a_1\cup a_2: a_1\in \idealI_1, a_2\in
\idealI_2\}$ is a maximal ideal of $\bool$, then there is an endomorphism
$h$ of $\bool$ such that $h\restriction \idealI_1={\id}_{{\mathcal
I}_1}$ and $h\restriction \idealI_2$ is constantly zero; 
but possibly there
are no such non--zero ideals $\idealI_1,\idealI_2$, 
 (then we call $\bool$
indecomposable).

In \S 2 we define the family of such endomorphisms (those defined by a 
schema and those defined by
a simple schema) and investigate this a little. Our main result (in \S 3) is
that for any $\lambda>\aleph_0$ there is a Boolean algebra of cardinality
$\lambda^{\aleph_0}$ (and even density character $\lambda$) with only endomorphisms
as above, of course there are $2^{\lambda^{\aleph_0}}$
such Boolean algebras
with no non trivial homomorphism from one to a distinct other (see \ref{3.1},
\ref{3.13}, \ref{3.14}); we also show that ``cardinality
$\lambda^{\aleph_0}$'' is a reasonable restriction 
(see \ref{3.15}, \ref{3.16AT}, \ref{4.20}). 

For simplicity, we concentrate on the case of $\cf(\lambda)>\aleph_{0}$;
note that this affect only the density character as 
$\cf(\lambda)=\aleph_0\Rightarrow \lambda^{\aleph_0}=
(\lambda^+)^{\aleph_0}$. 
How do we construct such $\bool$? The algebra $\bool$ extends the
Boolean algebra
$\bool_0$ 
which is 
freely generated by $\{x_\eta:\eta\in {}^{\omega{>}}\lambda\}$ and
is a subalgebra of its completion $\bool_0^c$. In fact, $\bool=\langle
\bool_0\cup
\{a_\alpha:\alpha<\alpha^*\}\rangle_{\bool_0^c}$, with $a_\alpha$ chosen
by induction on $\alpha$, has the form
$\bigcup\limits_n(d^\alpha_n\cap s^\alpha_n)$,
where $\langle d^\alpha_n:n<\omega\rangle$ is a maximal anti-chain of
$\bool$, for each $d_n^\alpha$ we have already decided that it will belong to $\bool$
and is based on ($=$ in the completion of the subalgebra generated by)
$\{x_\eta:\eta\in {}^{\omega{>}}\xi\}$ for some $\xi<\blackz(\alpha)<\lambda$,
and for some increasing $\eta_\alpha$ with the limit $\blackz(\alpha)$,
$s^\alpha_n\in \langle x_\nu:\eta_\alpha\restriction n\vartriangleleft\nu
\rangle_{\bool_0^c}$. Why these restrictions? We would like to ``kill''
undesirable endomorphisms and we shall omit appropriate countable
(quantifier free) types which the image of $a_\alpha$, if exists, has
to realize, so such restrictions give us tight control and so helps us
to ``diagonalize'' over all possible endomorphisms. To diagonalize we use
a black box --- it is presented in \S 1, but its existence is
not proved here (it is proved in \cite{Sh:309}).

\[*\qquad *\qquad *\]

In \cite{Sh:89}, answering a question of Monk, we have explicated the notion of
``a Boolean algebra with no endomorphisms except the ones induced by
ultrafilters on it'' (see \S 2 here) and proved the existence of one with
density character $\aleph_0$, assuming first $\diamondsuit_{\aleph_1}$ and
then only $CH$.  The idea was that if $h$ is an endomorphism of $\bool$, not
among the ``trivial'' ones, then there are pairwise disjoint $d_n\in \bool$
with
$h(d_n)\not\subseteq d_n$.  Then we can add, for some $S\subset \omega$, an
element $x$ such that $d_n\leq x$ for $n\in S, x\cap d_n=0$ for $n\not\in S$
while forbidding a solution for
\[\{y\cap h(d_n)=h(d_n):n\in S\}\cup \{y\cap h(d_n)=0:n\not\in S\}.\]
Later, further 
analysis had showed that the point is that we are 
omitting positive
quantifier free types. Continuing this, Monk succeeded to prove in ZFC, the
existence of such Boolean algebras of cardinality $2^{\aleph_0}$.  In his
proof he replaced some uses of the countable density character by the
$\aleph_1$--chain condition. Also, generally it is hard to omit
$<2^{\aleph_0}$ many types but because of the special character of the types
(as said above, 
positive quantifier free) and models involved, using $2^{\aleph_0}$ almost disjoint subsets
of $\omega$, he succeeded in doing this. Lastly, for another step in the
proof (ensuring idecomposability - see Definition \ref{2.1}) he (and
independently Nyikos) found it is in fact easier to do this when for every
countable set $Y\subseteq \bool$ there is $x\in \bool$ free over it.

The question of the existence of such Boolean algebras in other
cardinalities remained open (See \cite{DMR} and a preliminary list of problems
for the handbook of Boolean algebras by Monk).

In \cite{Sh:229} it is proved  (in ZFC) that there exist such $\bool$ of density character
$\lambda$ and cardinality $\lambda^{\aleph_0}$ 
whenever $\lambda>\aleph_0$; from this follows 
answers to some other 
questions from Monk's list, (combining
\ref{3.1} with \ref{2.5}).

Almost all the  present work 
is a revised version of \cite{Sh:229} 
but 
\ref{3.15} - \ref{4.20} 
were added; 
here as in \cite{Sh:229}   \S 2 repeats \cite{Sh:89}. 

\section{A black box}

\begin{Explanation}
\label{1.2}
{\rm
We shall let $\bool_0$ be the Boolean algebra freely generated by $\{\eta:\eta
\in {}^{\omega{>}}\lambda\}$, $\bool^c_0$ its completion and we can interpret
$\bool^c_0$ as a subset of ${\cal M}= {\cal H}_{<\aleph_1}(\lambda)$ (each
$a\in \bool^c_0$ has the form $\bigcup\limits_{n<\omega} s_n$ where $s_n$
is a Boolean combination of members of ${}^{\omega{>}}\lambda$). As
the $\eta\in {}^{\omega{>}}\lambda$ may be over-used, we replace them
for this purpose by $x_\eta$ (for example below let $F\in \tau_0$ 
be a unary
function symbol, $x_\eta=F(\eta)$).

Our desired Boolean algebra $\bool$ will be a subalgebra
of the competition $\bool^c_0$ of $\bool_0$ hence it extend $\bool_0$.
For our diagonalization, i.e. the omitting type, we need the
following case (we shall use $\kappa=\aleph_0$). That is we need a
family of subalgebras with endomorphism, for each we add an element
and promise to omit the type of the supposed image. The family is
sparse enough so that we can do it 
(i.e. with the different promises 
not 
 hindering one another too much), 
but dense enough so that every
endomorphism of the Boolean Algebra we construct is approximated. See
more accurate explanation in \ref{1.4}.
}
\end{Explanation}

\begin{convention}
\label{1.2d}
We fix $ \kappa \ge  {\aleph_0} $ for this section.
\end{convention}

\begin{Definition}
\label{1.3}
\begin{description}
\item[(1)] Let $\tau_n$, for $n<\omega$, be fixed vocabularies (=
signatures), $|\tau_n|\leq\kappa$, $\tau_{n}\subseteq \tau_{n+1}$, (with each
predicate and function symbol finitary for simplicity). Let $P_n\in
\tau_{n+1} \setminus \tau_n$ be unary predicates. Let ${\mathcal M}=
({\mathcal H}_{<\kappa^+}(\lambda),\in)$.
\item[(2)] For $ n<\omega$ let $\blackF_n$ be the family of sets (or
sequences) of the form $\{(f_\ell,N_\ell):\ell\leq n\}$ satisfying:
\begin{description}
\item[(a)] $f_\ell:{}^{\ell\geq}\kappa\longrightarrow {}^{\ell\geq}\lambda$
is a {\em tree embedding}, i.e.,
\begin{description}
\item[(i)] $f_\ell$ is length preserving, that is, $\eta$, $f_\ell(\eta)$ have
the same length,
\item[(ii)] $f_\ell$ is order preserving (of $\triangleleft$), 
moreover, for $\eta,\nu\in {}^{\ell
\geq} \kappa$ we have $\eta \triangleleft\nu$ iff $f_\ell(\eta)\triangleleft
f_\ell(\nu)$,
\end{description}
\item[(b)] $f_{\ell+1}$ extends $f_\ell$ (when $\ell+1\leq n$),
\item[(c)] $N_\ell$ is an $\tau'_\ell$--model of \power\ $\leq\kappa$,
$|N_\ell|\subseteq |{\cal M}|$, where $\tau'_\ell\subseteq\tau_\ell$,
\item[(d)] $\tau'_{\ell+1}\cap\tau_\ell=\tau'_\ell$ and $N_{\ell+1}
\restriction \tau'_\ell$ extends $N_\ell$,
\item[(e)] if $P_m\in \tau'_{m+1}$, then $P^{N_\ell}_m=|N_m|$ when $m<\ell
\leq n$, and
\item[(f)]  if $x,y\in N_\ell$ then $\{x,y\}\in N_ {\ell} $ 
 and $\emptyset \in N_ {\ell} $. 
\item[(g)] $\rRang(f_\ell)\subseteq N_\ell$
\end{description}
\item[(3)] Let $\blackF_\omega$ be the family of pairs $(f,N)$ such that for
some \sq\ $\langle(f_\ell,N_\ell):\ell< \omega\rangle$ 
the following hold:
\begin{description}
\item[(i)] $\{(f_\ell,N_\ell):\ell\leq n\}$ belongs to $\blackF_n$ for $n<
\omega$,
\item[(ii)] $f=\bigcup\limits_{\ell<\omega} f_\ell$, $N=\bigcup\limits_{n<
\omega} N_n$, (i.e., $|N|=\bigcup\limits_{n<\omega}|N_n|$, $\tau(N)=
\bigcup\limits_n \tau(N_n)$, and $N\restriction \tau(N_n)=\bigcup\limits_{n
<m<\omega}N_m\restriction \tau(N_n)$).
\end{description}
\item[(4)] For any $(f,N)\in\blackF_\omega$ let
 $\langle(f_n,N_n):n<\omega\rangle$ 
be as above (if $P_n 
\in \tau'_{n+1}$ for $n<\omega$ \then\ 
it
is easy to show that $(f_n,N_n)$ is uniquely determined by $(f,N)$- notice clauses
(d), (e) in (2)), so 
for each 
$(f^\alpha,N^\alpha)$ as in \ref{1.7A}
below 
$(f^\alpha_n,N^\alpha_n)$ for $n<\omega$
are defined as above. 
\item[(5)] A branch of $\Rang(f)$ or of $f$ (for $f$ as in (3)) is just any
$\eta\in {}^\omega\lambda$ such that for every $n<\omega$ we have $\eta
\restriction n\in\Rang(f)$.
\end{description}
\end{Definition}

\begin{Explanation of our Intended Plan}
\label{1.4}
{\rm
(of Constructing for example the Boolean algebra)

We will be given $\blackw=\{(f^\alpha,N^\alpha):\alpha<\alpha^* \}$,
so that every
branch $\eta $ of $f^\alpha$ converges to some $\blackz(\alpha)$, $\blackz(
\alpha)$ non-decreasing 
(in $\alpha$). We have a free object generated by
${}^{\omega{>}}\lambda$ (i.e., by $\langle x_\eta:\eta\in
{}^{\omega{>}}\lambda\rangle$, this is  $\bool_0$ in our case), and by
induction on $\alpha$ we define $\bool_\alpha$ and $a_\alpha$ for
$\alpha<\alpha^\ast$, such that $\bool_\alpha$ is increasing and continuous,
$\bool_{\alpha+1}$ is an extension of $\bool_\alpha$,
$a_\alpha\in \bool_{\alpha+1}
\setminus \bool_\alpha$ (usually $\bool_{\alpha+1}$ is generated by
$\bool_\alpha$ and
$a_\alpha$, and is included in the completion of $\bool_0$). Every element
will depend on few (usually $\leq \kappa$) members of
${}^{\omega{>}}\lambda$, and $a_\alpha$ ``depends'' in a peculiar way:
the set $Y_\alpha\subseteq{}^{\omega{>}} \lambda$ on which it
``depends'' is $Y^0_\alpha\cup Y^1_\alpha$, where $Y^0_\alpha$ is
bounded below $\blackz(\alpha)$ (i.e., $Y^0_\alpha\subseteq
{}^{\omega>}\zeta$ for some $\zeta<\blackz(\alpha))$ and $Y^1_\alpha$ is a
branch of $f^\alpha$ or something similar. See more in \ref{1.8}.
}
\end{Explanation of our Intended Plan}

\begin{Definition of the Game}
\label{1.5}
We define for $\blackw\subseteq\blackF_\omega$ a game $\Gm(\blackw)$,
which lasts $\omega$-moves.

In the $n$-th move:

Player II: Chooses $f_n$, a tree-embedding of $^{n\geq}\kappa$ into
${}^{n\geq}\lambda$, extending $\bigcup\limits_{\ell<n} f_\ell$, such that
$\Rang(f_n)\setminus\bigcup\limits_{\ell<n} \Rang(f_\ell)$ is disjoint to
$\bigcup\limits_{\ell<n}|N_\ell|$;

then

Player I: chooses $N_n$ such that $\{(f_\ell,N_\ell):\ell\leq n\}\in 
\blackF_n$.

In the end player II wins i \when 
$(\bigcup\limits_{n<\omega}f_n,
\bigcup\limits_{n<\omega}N_n)\in \blackw$.
\end{Definition of the Game}

\begin{Remark}
\label{1.6}
{\rm
We shall be interested in $\blackw$ such that player II wins (or at least does not
lose) the game, but $\blackw$ is ``thin''. Sometimes we need a strengthening of
the first player in two respects: he can demand (in the $n$-th move)
that $\Rang(f_{n+1})\setminus\Rang(f_n)$ is outside a ``small'' set, and in
the zero move he can determine an arbitrary initial segment of the play.
}
\end{Remark}

\begin{Definition}
\label{1.7}
We define, for $\blackw\subseteq \blackF_\omega$, a game
$\Gm'(\blackw)$ which lasts $\omega$-moves.

In the zero move:

\noindent Player I  chooses $k<\omega$ and $\{(f_\ell,N_\ell):\ell\leq k\}
\in \blackF_k$, and $X_0\subseteq {}^{\omega{>}}
\lambda$, $|X_0|<\lambda$.

In the $n$-th move, $n>0$:

\noindent Player II chooses $f_{k+n}$, a tree embedding of ${}^{(k+n)\geq}
\kappa$ into ${}^{(k+n)\geq}\lambda$, with $\Rang(f_{k+n})\setminus
\bigcup\limits_{\ell<k+n} \Rang(f_\ell)$ disjoint to $\bigcup\limits_{\ell<
k+n} N_\ell\cup\bigcup\limits_{\ell<n} X_\ell$.

\noindent Player I chooses $N_{k+n}$ such that $\{(f_\ell,N_\ell):\ell \leq
k+n\}\in \blackF_{k+n}$ and he chooses $X_n\subseteq {}^{\omega{>}}\lambda$ satisfying
$|X_n|<\lambda$.

In the end of the play, player II wins 
\when\ 
$(\bigcup\limits_{n<\omega}f_n,
\bigcup\limits_{n<\omega} N_n)\in \blackw$ 
\end{Definition}

\begin{Remark}
\label{1.8}
{\rm
What do we want from $\blackw$? First that by adding an element
(to $\bool_0$) for each $(f,N)$, we can ``kill'' every undesirable
endomorphism, for this $\blackw$
has to 
``encounter" 
every possible endomorphism, and this will be served by
``$\blackw$ a barrier'' defined below. For this $\blackw=\blackF_\omega$ is O.K.
but we also want $\blackw$
to be thin enough so that various demands will have small interaction. For
this, disjointness and more are demanded.
}
\end{Remark}
\medskip

\begin{Definition}
\label{1.6A}
\begin{description}
\item[(1)] We call $\blackw\subseteq \blackF_\omega$ a {\em strong barrier} \uif\
Player II wins in the game $\Gm(\blackw)$ and even $\Gm'(\blackw)$
(which just means he has a winning strategy.)
\item[(2)] We call $\blackw$ a {\em barrier} \uif\ Player I
does not win in the game
$\Gm(\blackw)$ and even does not win in $\Gm'(\blackw)$.
\item[(3)] We call $\blackw$ disjoint \uif\ for any distinct
$(f^\ell,N^\ell)\in \blackw$
($\ell=1,2$), $f^1$ and $f^2$ have no common branch.
\end{description}
\end{Definition}

\begin{The Existence Theorem}
\label{1.7A}
\begin{description}
\item[(1)] If $\lambda^{\aleph_0}=\lambda^{\kappa}$, $\cf(\lambda)>\aleph_0$
\then\ there is a strong disjoint barrier.
\item[(2)] Suppose $\lambda^{\aleph_0}=\lambda^{\kappa}$, $\cf(\lambda)>
\aleph_0$. \Then\ there is
\[\blackw=\{(f^\alpha,N^\alpha):\alpha<\alpha^*\}\subseteq 
\blackF_\omega\]
and a non-decreasing function $\blackz:\alpha^*\longrightarrow\lambda$ such
that:
\begin{description}
\item[(a)] $\blackw$ is a strong disjoint barrier, moreover for every stationary
$S\subseteq\{\delta<\lambda:\cf(\delta)=\aleph_0\}$, the set
$\{(f^\alpha,N^\alpha):
\alpha<\alpha^*,\blackz (\alpha)\in S\}$ is a disjoint barrier,
\item[(b)] $\cf(\blackz(\alpha))=\aleph_0$ for $\alpha<\alpha^*$,
\item[(c)] every branch of $f^\alpha$ is an increasing sequence converging to
$\blackz(\alpha)$,
\item[(d)] each $N^\alpha_n$ is transitive, i.e.: if ${\mathcal
M}\models$``$a\in b$'', $b\in N^\alpha_n$, $b \notin \lambda$,
 \then\ $a\in N^\alpha_n$,
(we call $\{(f_\ell,N_\ell):\ell\leq n\}$, transitive if each $N_\ell$ is 
transitive and similarly $\{(f_\ell,N_\ell):\ell<\omega\}$ and $\blackW$).
\item[(e)] if $\blackz(\beta)=\blackz(\alpha)$, $\beta+\kappa^{\aleph_0}\leq
\alpha<\alpha^*$ and $\eta$ is a branch of $f^\alpha$, then $\eta
\restriction k\not\in N^\beta$ for some $k<\omega$,
\item[(f)] when $\lambda=\lambda^\kappa$ we can demand: if $\eta$ is a branch
of $f^\alpha$ and $\eta\restriction k\in N^\beta$ for all $k<\omega$ (where
$\alpha,\beta<\alpha^*$) then $N^\alpha\subseteq N^\beta$ (and even for
every $n<\omega$, $N^\alpha_n\in N^\beta$).
\end{description}
\end{description}
\end{The Existence Theorem}

Proof: See  
     \cite[3.11]{Sh:309}, 
     \cite[3.16]{Sh:309}.

\section{Preliminaries on Boolean Algebras}
\medskip
We present here some easy material.

\begin{Definition}
\label{2.1}
\begin{description}
\item[(1)] For any endomorphism $h$ of a Boolean Algebra $\bool$, let
\[\hspace{-0.7cm}\begin{array}{lcl}
{\ExKer}(h)&=&\{x_1\cup x_2:\ h(x_1)=0\mbox{, and } h(y)=y\mbox{ for every }y\leq x_2\},\\
{\ExKer}^*(h)&=&\{x\in \bool:\mbox{ in }\bool/
\ExKer(h)\mbox{, below }x/\ExKer
(h),\\
&&\qquad\qquad\mbox{there are only finitely many elements}\}.
  \end{array}\]
\item[(2)] A Boolean algebra is endo-rigid \uif\  for every
endomorphism $h$ of $\bool$, $\bool/\ExKer(h)$ is finite (equivalently:
$1_\bool\in \ExKer^\ast (h)$).
\item[(3)] A Boolean algebra is indecomposable \uif\  there are
no two disjoint ideals ${\mathcal I}_0$, ${\mathcal I}_1$ of $\bool$ 
(except $0_\bool$ of course), each
with no maximal member, which generate a maximal ideal 
of $\bool$ (that is $\{a_0\cup
a_1:\,a_0\in {\mathcal I}_0, a_1\in {\mathcal I}_1\}$).
\item[(4)] A Boolean algebra $\bool$ is $\aleph_1$-compact 
\uif\  
for every pairwise disjoint $d_n\in \bool$ (for $n<\omega$)
for some $x\in \bool$, we have
$x\cap d_{2n+1}=0$, $x\cap d_{2n}=d_{2n}$.
\end{description}
\end{Definition}

\begin{Lemma}
\label{2.2}
\begin{description}
\item[(1)] A Boolean algebra $\bool$ is endo-rigid \uiff\ every
endomorphism of $\bool$ is the endomorphism of some scheme (see Definition
\ref{2.3}(1),(3) below).
\item[(2)] A Boolean algebra $\bool$ is endo-rigid and indecomposable
\uiff\ every endomorphism of 
$\bool$ is the endomorphism of some
simple scheme (see Def \ref{2.3}(2) below).
\item[(3)] For every scheme of an endomorphism of $\bool$ 
\there\
is one and only one endomorphism of the scheme.
\end{description}
\end{Lemma}

\par\noindent
Proof. \quad Easy.

\begin{Remark}
\label{2.2A}
{\rm
\begin{enumerate}
\item In fact, for a Boolean algebra $\bool$, we have $\{h: h$ is an
endomorphism of $\bool$ defined by a scheme$\}$ is a sub-semi-group of ${\rm
End}(\bool)$, even a normal one (as $(\bool,{\rm End}(\bool))$
is interpretable in ${\rm End}(\bool)$).
\item Similarly for simple schemes.
\end{enumerate}
}
\end{Remark}

\begin{Definition}
\label{2.3}
\begin{description}
\item[(1)] A scheme of an endomorphism of $\bool$ consists of a partition $a_0,
a_1,b_0,\ldots,b_{n-1}$, $c_0,\dots,c_{m-1}$ in $\bool$ of $1_\bool$,
with maximal non-principal ideals $\idealI_\ell$ below $b_\ell$
for $\ell<n$ (in other words $\idealI_\ell$ is a maximal ideal of
$\bool \rest b_\ell$) and  non-principal
 ideals $\idealI^0_\ell,\idealI^1_\ell$ below $c_\ell$ for $\ell<m$ such that
$\idealI^0_\ell \cup \idealI^1_\ell$ generates a
maximal non principal ideal below
$c_\ell$ and $\idealI^0_\ell \cap \idealI^1_\ell = \{0_\bool\}$,
a number $k\leq n$,
and a partition $b^*_0,\ldots,b^*_{n-1}, c^*_0,\ldots,c^*_{m-1}$ of
$a_0\cup b_0\cup \dots\cup b_{k-1}$. We assume also that
\[[k+m>0\ \ \Rightarrow\ \ a_0=0],\qquad [(n-k)+m>0\ \ \Rightarrow\ \
a_1=0]\]
and except possibly $a_0,a_1$ there are no zero elements in the
partition $a_0, a_1$, $ b_0, \ldots, b_{n-1}, c_0, \ldots, c_{m-1}$.
\item[(2)] The scheme is simple if $m=0$.
\item[(3)] The endomorphism of the scheme is the unique endomorphism $h:\bool
\longrightarrow \bool$ such that:
\begin{description}
\item[(i)]   $h(x)=0$ when $x\leq a_0$ or $x\in \idealI_\ell$, $\ell<k$,
or $x\in \idealI^0_\ell$, $\ell<m$,
\item[(ii)]  $h(x)=x$ when $x\leq a_1$ or $x\in \idealI_\ell$, $k\leq\ell<n$ or
$x\in \idealI^1_\ell$, $\ell<m$,
\item[(iii)] $h(b_\ell)=b_\ell^*$ when $\ell<k$,
\item[(iv)]  $h(b_\ell)=b_\ell\cup b^*_\ell$ when $k\leq \ell<n$,
\item[(v)]   $h(c_\ell)=c_\ell\cup c^*_\ell$ when $\ell<m$.
\end{description}
\end{description}
\end{Definition}

So, an endomorphism of a scheme is a ``trivial" endomorphism defined by
ideals, essentially maximal ones, and finitely many elements.

\begin{Claim}
\label{2.4}
\begin{enumerate}
\item If $h$ is an endomorphism of a Boolean Algebra
$\bool$, and $\bool/ \ExKer(h)$
is infinite \then\  there are pairwise disjoint
$d_n\in \bool$ (for $n<\omega$)
such that $h(d_n)\nleq d_n$.
\item We can demand that: $h(d_n)\cap d_{n+1}\neq 0$, and if $\bool$ satisfies
the c.c.c.,\then\ $\{d_n:\,n<\omega\}$ is a maximal antichain.
\end{enumerate}
\end{Claim}

\par\noindent
Proof. (1)\quad As $\bool/\ExKer(h)$ is infinite we can choose inductively
$d_n\in \bool$ such that $d_n\notin \ExKer(h)$, $[\ell<n\ \Rightarrow\
d_\ell\cap d_n=0_\bool]$ and $\{x/\ExKer(h): x\in \bool\ \&\
x\cap\bigcup\limits_{\ell\leq n} d_\ell=0_\bool\}$ is infinite.
It is enough for
each $n$ to find $d^*_n\leq d_n$ 
such that $h(d^*_n)\not\leq d^*_n$. Since
$d_n\notin \ExKer(h)$, clearly (by the definition of $\ExKer(h)$)
we have $h(d_n)>0_\bool$ and for some $d_n'\leq d_n$, $h(d_n')\neq d_n'$.
\par\noindent
 \underline{Case 1:}  $h(d_n)\not\leq d_n$, let $d^*_n=d_n$.
\par \noindent
\underline{Case 2:}
$h(d_n)=d_n$. 

Now if $h(d_n')\not\leq d_n'$ let $d^*_n=d_n'$ and otherwise 
$h(d'_n)\leq d'_n$ so by the choice of $d'_n$ we have 
$h(d_n')<h(d_n)$,
let $d^*_n=d_n-d_n'$, so $h(d^*_n)=h(d_n)-h(d_n')=d_n-h(d_n')>d_n-d_n'
=d^*_n$ so $d^*_n$ is as
required.
\par \noindent \underline{Case 3:}
Neither case 1 nor case 2. 

So $h(d_n)\leq d_n$ but $h(d_n)\neq d_n$
hence $h(d_n)<d_n$. So $h(d_n- h(d_n))\leq h(d_n-0_\bool)=h(d_n)$ is
disjoint from $d_n-h(d_n)$, so if $h(d_n-h(d_n))>0_\bool$ we let
$d_n^*=d_n-h(d_n)$. So assume not, so
$d_n-h(d_n)\in\Ker(h)\subseteq \ExKer(h)$, and hence $h(h(d_n))=
h(0_\bool \cup h(d_n))=h((d_n-
h(d_n))\cup h(d_n))=h(d_n)$ and necessarily $h(d_n)\notin \ExKer(h)$
(as $d_n\notin \ExKer(h)$, $d_n - h(d_n)\in \ExKer(h)$), hence case 2
apply to $h(d_n)$ and we are done.

\noindent (2)\quad Let $c_n=h(d_n)-d_n>0$, so $m\neq n\ 
\Rightarrow\ d_n\cap d_m=0_\bool \Rightarrow c_m\cap
c_n=0_\bool$ and
\[ d_n\cap c_n= d_n\cap (h(d_n)-d_n) 
\leq d_n \cap (-d_n) = 0_\bool\]
so $d_n\cap c_n=0_\bool$.

By Ramsey theorem,\wolog, for all $m<n$ the truth
value of $d_m\cap c_n=0_\bool$ is the same and of 
$c_n\cap d_m=0_\bool$ is the same.

\par \noindent
Now we prove
\begin{description}
\item[$(*)$] for some $\langle d'_n: n<\omega\rangle$ 
we have $d'_n\in  
\bool$, $h(d'_n)\not\leq d'_n$ moreover $h(d_n')\cap d_{n+1}'>0$
and $n<m \Rightarrow d'_n \cap d'_m= 0_\bool$.
\end{description}

\par \noindent
\underline{Case 1:} $c_0\cap d_1 > 0_\bool$.

Let $d_n'=d_{2n+2}\cup (c_{2n}\cap d_{2n+1})$; now $\langle d_n':n<\omega
\rangle$ are pairwise disjoint as the $d_n$'s are. Now as $h(d_m) \geq
c_m$ for $m<\omega$ clearly
\[h(d_n')\geq h(d_{2n+2})\geq c_{2n+2}\geq c_{2n+2}\cap d_{2n+3}>0_\bool\]
$$
d_{n+1}'\geq c_{2n+2}\cap d_{2n+3}>0_\bool,
$$
so\quad  $h(d_n')\cap d_{n+1}'\geq c_{2n+2}\cap d_{2n+3}>0_\bool$.
\par \noindent
\underline{Case 2:}
$c_1\cap d_0>0_\bool$.

Let $d_n'=d_{2n+3}\cup (c_{2n+1}\cap d_{2n})$. Now $\langle d_n':
n<\omega\rangle$ are pairwise disjoint (as $\langle d_{2n}\cup
d_{2n+3}:n<\omega\rangle$ are), $h(d'_n)\geq h(d_{2n+3})\geq
c_{2n+3}\geq c_{2n+3}\cap
d_{2n+2}>0$ and $d_{n+1}^\prime\geq c_{2(n+1)+1} \cap d_{2(n+1)}=
c_{2n+3}\cap d_{2n+2}>0$. So clearly $h(d_n')
\cap d_{n+1}^\prime\geq c_{2n+3}\cap d_{2n+2}>0$.

\par \noindent
\underline{Case 3:} Neither Case 1 nor Case 2.

As we have noted above $d_n\cap c_n=0_\bool$ by the case assumption's we
have $d_n\cap c_m=0_\bool$ for every $m,n<\omega$ 
and of course $n\neq m \Rightarrow d_n\cap d_m
= 0_\bool\ \&\ c_n\cap c_m=0_\bool$. Lastly let 
$d_n' = d_{n+1} \cup c_n$, they are as required 
e.g. $h(d'_n)\cap d'_{n+1}=(h (d_{n+1})\cup h(c_n))\cap
(d_{n+2}\cup c_{n+1})\geq h(d_{n+1})\cap c_{n+1}=c_{n+1}>0_\bool$.

So we have proved $(*)$. Now renaming $d_n'$ as $d_n$, $\langle d_n: n<
\omega\rangle$ satisfies (part (1) and) the first demand of part (2).

If $\bool$ satisfies the c.c.c., we can find $\alpha\in [\omega,\omega_1)$ and
$d_\beta$ for $\beta\in [\omega,\alpha)$ such that $\langle d_\beta:\beta<
\alpha\rangle$ is a maximal antichain of $\bool$, 
\wolog,
$\alpha\leq \omega+\omega$. Now let $d_n'$ be $d_n\cup d_{\omega+n}$ if
$\omega+n<\alpha$, and $d_n$ otherwise. So $n\neq m\ \Rightarrow\ d_n'\cap
d_m'=0$ and $h(d_n')\cap d_{n+1}'\geq h(d_n)\cap d_{n+1}>0$, so $\langle d_n':n<\omega\rangle$ is
as required. \qquad\hfill $\qed_{\ref{2.4}}$

\begin{Definition}
\label{2.4u}
A Boolean Algebra $\bool$ is Hopfian 
\uif\  
every onto endomorphism
of $\bool$ is one-to-one. A Boolean 
Algebra $\bool$ is dual Hopfian 
\uif\ 
every one to one endomorphism is onto.
\end{Definition}

\begin{Lemma}
\label{2.5}
\begin{description}
\item[(1)] Every atomless endo-rigid Boolean Algebra $\bool$ is Hopfian and dual
Hopfian.
\item[(2)] Also $\bool+\bool$ is Hopfian 
(and dual Hopfian), however it is  not rigid.
\end{description}
\end{Lemma}
\medskip

Proof: Easy to check using \ref{2.2}, \ref{2.3}.

\section{The Construction}

\begin{Main Theorem}
\label{3.1}
Suppose $\cf(\lambda)>\aleph_0$. \Then\ there is a 
Boolean algebra $\bool$ such that:
\begin{description}
\item[(a)] $\bool$ satisfies the \Ccc\ and is atomless,
\item[(b)] $\bool$ has power $\lambda^{\aleph_0}$ and 
has algebraic density $\lambda$ 
(in the 
Boolean cardinal 
invariant notation, $\pi(\bool)=\lambda$), 
 this means:  
\[\min\{|X|: X\subseteq \bool\setminus\{0 _ \bool 
\}
\mbox{ and } (\forall y\in \bool)
(\exists x\in X)(y>0\ \Rightarrow\ x\leq y)\},\]
\item[(c)] $\bool$ is endo-rigid and indecomposable.
\end{description}
\end{Main Theorem}
\medskip

\par\noindent
Proof. Let $\tau_n$ 
for $ n < \omega $ 
be as in \S 1 for $\kappa = \aleph_0$,
we use $ \tau'_n $  
 with
$\tau'_0$ 
 having unary predicate $Q$, binary predicate $\leq$,
individual constants 1, 0, binary function symbols $\cup,\cap,-$ and unary
function symbol $H$ (and more)
and 
$ P_n \notin \tau'_{n+1 } \setminus \tau'_n $.  
 We shall use Theorem \ref{1.7A}(2) for $\lambda$
and $\kappa=\aleph_0$, and let $\blackw=\{(f^\alpha,N^\alpha):\,
\alpha<\alpha^*
\}$, the function $\blackz$, the model ${\cal M}=({\cal H}_{<\aleph_1}
(\lambda),\in)$ and ${\mathcal T}={}^{\omega>}
\lambda$ be as there. We call $\alpha<\alpha^*$ a candidate if
\[\bool^{N^\alpha}=\bool[N^\alpha]=(Q^{N^\alpha},1^{N^\alpha},0^{N^\alpha},
\cup^{N^\alpha}, \cap^{N^\alpha},-^{N^\alpha},\leq^{N^\alpha})\]
is a Boolean algebra and $h_\alpha=H^{N^\alpha}\restriction \bool^{N^\alpha}$
is an endomorphism of $\bool^{N^\alpha}$; 
of course $\cup^{N^\alpha}$ means
$\cup^{N^\alpha}\restriction Q^{N^\alpha}$ etc and we are 
demanding that all the relevant predicates and function symbols belongs 
to $\tau_{N^\alpha}$. 

We will think of the game as follows: Player I tries to produce a non
trivial endomorphism $h$. Player II supplies (via the range of $f_\ell$)
elements in $\bool_0$ (see Stage A below) 
and challenges Player I for defining $h$ on them. So
Player I plays models $N_\ell$ 
in the vocabulary $\tau'_\ell$ 
which is mainly a subalgebra of the Boolean 
algebra we are constructing, with additional 
elements and expanded  
by, 
in particular, the  
distinguished function symbol  
$ H \in \tau'_0 $ 
which is interpreted as
an endomorphism of Boolean
algebras. In the end, as $\blackw$ 
is a  
barrier,  for some such play 
 we will get a model
$N^\alpha\in \blackw$, in the 
vocabulary 
$\bigcup\limits_{\ell<\omega}\tau'_\ell$
which includes a function symbol $H$. 
We can think of $N^\alpha$ as a Boolean
algebra $\subseteq \bool^c_0$ with an endomorphism $h_\alpha=H^{N_\alpha}$.
\medskip

\noindent\underline{Stage A}\quad Let $\bool_0$ be the Boolean algebra freely
generated by $\{x_\eta:\, \eta\in {}^{\omega{>}}\lambda\}$, and $\bool^c_0$ be
its completion. For $A\subseteq \bool^c_0$ let
$\langle A\rangle_{\bool^c_0}$ be the
Boolean subalgebra of $\bool_0^c$ that $A$ generates. As $\bool_0$
satisfies the
c.c.c. every element of $\bool^c_0$ can be represented as a countable union of
members of $\bool_0$, and as $\bool_0$ is free we get  
$\|\bool_0^c\|=\lambda^{\aleph_0}$. We say $x\in \bool^c_0$
is based on (or
supported by) $J\subseteq {}^{\omega>}\lambda$ if it is based on (or
supported by) $\{x_\nu:\,\nu\in J\}$ that is $\bool^c_0\models \mbox{``}
x=\bigcup\limits_{n<\omega} y_n\mbox{''}$, where each $y_n$ is in the
subalgebra generated by $\{x_\nu:\, \nu\in
J\}$; we shall also say that $J$ is a support of $x$. Let ${\rm supp}(x)$
be the minimal such $J$; it is easy to prove its existence.
[Why? Let $x=\bigcup\limits_{n<\omega} y_n$, where $y_n=
\sigma_n(\ldots, x_{\eta_{n, \ell}}'\ldots)_{\ell\leq k_n}$; as if
$y_n=\bigcup\limits_{\ell<k} y_{n, \ell}$ we can replace $y_n$ by
$y_{n, 0}, \ldots, y_{n, k-1}$,
hence 
 \wolog, for each
$n$, for some disjoint finite $u_n, v_n\subseteq {}^{\omega>}\lambda$ we
have $y_n = \bigcap\limits_{\eta\in u_n} x_\eta \cap
\bigcap\limits_{\eta\in v_n} (-x_\eta)$. Also we can replace $u_n$ by
any $u\subseteq u_n$ such that $y'=\bigcap\limits_{\eta\in u}
x_\eta \cap \bigcap\limits_{\eta\in v_n} (-x_\eta)$ satisfies $y'\leq x$.
So \wolog
$$
u\subseteq u_\eta\ \&\ u\neq u_\eta
\Rightarrow \bigcap\limits_{\eta\in u} x_\eta \cap
\bigcap\limits_{\eta\in v _n } 
 (-x_\eta) - x >0.
$$
Similarly \wolog
$$
v\subseteq v_n\ \&\ v\neq v_n\ \Rightarrow\ \bigcap\limits_{\eta\in
u_n} x_\eta \cap \bigcap\limits_{\eta\in v} (-x_\eta)-x>0.
$$
Lastly let $J= \bigcup\limits_{n<\omega} u_n \cup\bigcup\limits_{n<\omega} v_n$,
clearly $J$ is a support of $x$. If $J$ is not minimal then let $J'$
be a support of $x$ \st\ $J\not\subseteq J'$ as witness e.g. by
$\langle y'_n:n<\omega\rangle$. So  for some $n$, $u_n\cup
v_n\not\subseteq J'$, by symmetry \wolog\ 
$u_n\not\subseteq J'$, but then
$u'=u_n\cap J'$ contradicts the statement above.]
\par \noindent
\Wlog, not only $\bool_0^c\subseteq {\mathcal M}$ but
$x\in \bool_0^c$ implies
that the transitive closure of $\{x\}$ in ${\mathcal M}$ includes 
${\rm supp}(x)$. 
We shall now define by induction on $\alpha<\alpha^*$, the truth
value of ``$\alpha\in Y_\ell$'' ($\ell=1,2,3$), ``$\alpha\in Y$'',
``$\alpha\in Y'$'', the sequence
$\eta_\alpha$, and members $a_\alpha$,
$b^\alpha_m$, $c^\alpha_m$, $d^\alpha_m$, $s^\alpha_m$ of $\bool^c_0$ for
$m<\omega$ when $\alpha\in Y'$ such that,
$Y\cup Y_1\cup Y_2\cup Y_3 \subseteq Y'$ and
letting $\bool_\alpha=\langle \bool_0\cup\{a_\gamma: \gamma<\alpha,\;
\gamma\in Y'\}\rangle_{\bool^c_0}$ 
we have \footnote{ actually 
the $\alpha\in Y'\setminus Y= Y'\setminus Y_1$ 
have no real role here, but have in \ref{3.13} later.}:

\begin{description}
\item[$\circledast_\alpha$]

\begin{enumerate}
\item[(1)] $\eta_\alpha$ is a branch of 
$\Rang(f^\alpha)$, and  $\eta_\alpha\neq
\eta_\beta$ for $\beta<\alpha$;
\item[(2)] if $\alpha\in Y'$, then for some $\xi<\blackz(\alpha)$:

$a_\alpha=\bigcup\limits_{m<\omega}(s^\alpha_m\cap d^\alpha_m)$, where
$\langle d^\alpha_m:\,m<\omega\rangle$ is a maximal antichain of non zero
elements (of $\bool^c_0$), $d^\alpha_n\in \bool_\alpha$, $\bigcup\limits_m{\rm
supp}(d^\alpha_m)\subseteq
{}^{\omega>}\xi$, $s^\alpha_m\in \langle x_\rho:\ \eta_\alpha\restriction m
\vartriangleleft \rho,\ \rho\in {}^{\omega{>}}\lambda\rangle_{\bool^c_0}$, and
$d^\alpha_m>s^\alpha_m\cap d^\alpha_m> 0$, and $b^\alpha_n, c^\alpha_n\in
\bool_\alpha$ are based on ${}^{\omega{>}}\blackz(\alpha)$ and
$a_\alpha\notin \bool_\alpha$;
\item[(3)] if $\alpha\in Y$, then $b^\alpha_n$, $d^\alpha_n\in N^\alpha_0$,
$c^\alpha_n$, $s^\alpha_n\in N^\alpha$ (hence, by clause $(g)$ of 
Definition \ref{1.3}(2) each is based on $\{x_\nu:\,
\nu\in {}^{\omega{>}}\lambda,\ \nu\in N^\alpha\}$), and $b^\alpha_n\cap
b^\alpha_m=0$ for $n\neq m$;
\item[(4)] if $\beta<\alpha$, $\beta\in Y$, then $(\beta\in Y'$ and)
$\bool_\alpha$ omits the type
$$
p_\beta=\{x\cap b^\beta_n=c_n^\beta:\ n<\omega\}.
$$
\end{enumerate}
\end{description}
\medskip

Before we carry out the construction observe:
\medskip

\begin{Crucial Fact}
\label{3.2}
For any $x\in \bool_\alpha$ letting $\zeta=\blackz(\alpha)$
there are a finite subset $J$ of
${}^{\omega>}\lambda$, $k<\omega$, $\xi<\zeta$, and $\alpha_0<
\ldots<\alpha_{k-1}<\alpha$ such that
\begin{description}
\item[(a)] $\blackz(\alpha_0)=\blackz(\alpha_1)=\blackz(\alpha_2)=\dots=
\blackz(\alpha_{k-1})=\zeta$,
\item[(b)] $x$ is based on
\[\{x_\nu:\, \nu\in J\cup {}^{\omega>}\xi\mbox{ or }\nu\in {\rm supp}(
s^{\alpha_\ell}_m)\mbox{ for some }\ell< k,\ m<\omega\}.\]
\item[(c)] $x=\sigma(a_{\alpha_0},\ldots,a_{\alpha_{k-1}}, 
b_0,\ldots, b_{n-1})$, 
for some Boolean term $\sigma$,\\ and $b_0, \ldots,
b_{n-1}\in \langle \bool_0\cup \{a_\alpha: 
\blackz(\alpha)< \xi\}\rangle$,
\\and if $x\in \bool_0$ then $k=0$ and
$n$ is minimal.
\end{description}
\end{Crucial Fact}
\medskip

\noindent \underline{Continuation of the proof of \ref{3.1}}
\noindent \underline{Stage B}\quad Let us carry out the
construction on $\alpha$. For
$\xi<\lambda$, $w\subseteq\alpha^*$ let
\[I_{\xi,w}=\{\nu:\, \nu\in {}^{\omega>}\xi\mbox{ or }\ \nu\in
\bigcup_{{m<\omega},{\gamma\in w}}{\rm supp}(s^{\gamma}_m)\}.\]
We call $\alpha$ a good candidate if $(\alpha<\alpha^*$ and)
$\bool[N^\alpha]$ is a subalgebra of $\bool_\alpha$, $x\in 
\bool [N^\alpha]\Rightarrow \supp (x) \subseteq N^\alpha$ of course 
 and $h_\alpha=H^{N_\alpha}\rest \bool^{N_\alpha}$ is an
endomorphism of $\bool[N^\alpha]$ (note that $h_\alpha$ maps
$\bool [N^\alpha_n]$ into $\bool [N^\alpha_n]$ for $n<\omega$). We let $\alpha\in Y_1$ 
\uiff
\begin{enumerate}
\item[$\otimes^1_\alpha$]
\begin{description} 
\item[$(\alpha)$] $\alpha$ is a good candidate 
\item[$(\beta)$] there are $d^\alpha_m\in
N^\alpha_0\cap \bool_\alpha$ for $m<\omega$, $d^\alpha_m\neq 0$,
$d^\alpha_m\cap d^\alpha_\ell=0$ for $m\neq \ell$, such that
$\langle d^\alpha_m:m<\omega\rangle$ is a maximal antichain of
$\bool^c_0$ and for some $\xi<\blackz(\alpha)$
each $d^\alpha_m$ is based on ${}^{\omega>}\xi$, and there are a branch
$\eta_\alpha$ of $\Rang(f^\alpha)$ and $s^\alpha_m\in N^\alpha\cap 
\bool_\alpha$ ($m<\omega$)
as in (1), (2) above,
\item[$(\gamma)$] in addition if we add $\bigcup\limits_{n<\omega}(
s^\alpha_n\cap d^\alpha_n)$ to $\bool_\alpha$ then
\begin{description}
\item[\qquad (a)]  each $p_\beta$ ($\beta\in Y\cap\alpha$)
is still omitted
\item[\qquad (b)$_1$]
$p_\alpha=:\{x\cap h_\alpha(d^\alpha_m)=
h_\alpha(d^\alpha_m\cap s^\alpha_m):
m<\omega\}$ is omitted.
\end{description}
\end{description}
\end{enumerate}
${}$

Let $\alpha\in Y_2$ \uiff\  $\alpha\notin Y_1$ and
$\otimes^2_\alpha$ holds where
\begin{enumerate}
\item[$\otimes^2_\alpha$] is defined like
 $\otimes^1_\alpha$ replacing clause $({\bf b})_1$ by
\begin{description}
\item[\qquad\quad  (b)$_2$] the type $p'_\alpha$ is realized where
$p'_\alpha=\{x\cap h_\alpha(d_{2n}\cup d_{2n+1})= h_\alpha(d_{2n}):
n<\omega\}$ and $h_\alpha(s_{2n})=1$, $h_\alpha(s_{2n+1})=0$ and $\langle
d_n: n<\omega\rangle \in N^\alpha_0$ but $\bool[N^\alpha]$ omit $p'_\alpha$.
\end{description}

Let $\alpha\in Y_3$ iff $\alpha\notin Y_1\cup Y_2$ and
\item[$\otimes^3_\alpha$] we have $(\alpha)$ +
$(\beta)$ + $(\gamma)$(a) from $\otimes^1_\alpha$.
\end{enumerate}

Let $Y=Y_1$, $Y'=Y_1\cup Y_2\cup Y_3$ and for $\alpha\in Y_\ell$ let
$\otimes_\alpha$ mean $\otimes^\ell_\alpha$.

If $\alpha\in Y'$ we choose $\eta^\alpha$, $d^\alpha_n$, $s^\alpha_m$,
satisfying $\otimes_\alpha$ so also $\bool_{\alpha+1}$ is well defined 
and if $\ell=1$ let $b^\alpha_m=h_\alpha(d^\alpha_m)$,
$c^\alpha_m =h_\alpha(d^\alpha_m\cap s^\alpha_m)$ for $m<\omega$ 
, if $\ell\in\{2,3\}$ we can still choose $b_m^\alpha, c^\alpha_m$ for
$m<\omega$ \st\ $\circledast_\alpha$ holds (e.g. $\{\langle n\rangle:n<\omega\}
\subseteq N^\alpha_0$ by clause (f) of \ref{1.3}(2), so $b^\alpha_n=
x_{\langle 2n+1\rangle}-\bigcup\limits_{m<n} x_{\langle 2n+1\rangle}, 
c^\alpha_n=b^*_n\cap x_{\langle 2n\rangle}$). 

If $\alpha\in \alpha^*\setminus Y'$ we leave $a_\alpha$, $\eta_\alpha$ and $d^\alpha_n$,
$s^\alpha_n$ (for $n<\omega$) undefined, and so $\bool_{\alpha+1}
= \bool_\alpha$. So we have carried the induction.

So ``$\alpha\in Y'$'' means that Player I played Boolean Algebras and
endomorphisms as in the previous remark and we get in the end a Boolean
Algebra with the same properties.

The desired Boolean algebra $\bool$ is $\bool_{\alpha^*}=
\cup\{\bool_\alpha:\alpha<\alpha^*\}$. We shall investigate it
and eventually prove that it is endo-rigid (in \ref{3.11}) 
and indecomposable  (in \ref{3.12}) thus proving $\bool_{\alpha^*}$ is as
required in clause (c) of  \ref{3.1}, where (\ref{3.1}(a),
\ref{3.1}(b) hold trivially noting that $|\bool_{\alpha^*}|$
is $\leq|\bool_0^c|\leq |\bool_0|^{\aleph_0}=\lambda^{\aleph_0}$
and is $\geq |Y'|\geq \lambda^{\aleph_0}$ which will be proved later (see \ref{3.12AT})
and $a_\alpha \notin \bool_\alpha$ by (2) from stage A). The rest of the
proof is broken to facts and claims in this framework.
\medskip

Note also
\begin{Fact}
\label{3.3}
\begin{description}
\item[(1)] For $\nu\in {}^{\omega>}\lambda$, 
$x_\nu$ is free over $\{x_\eta:\
\eta \in{}^{\omega>}\lambda,\eta \neq \nu\}$ in $\bool_0$  
hence also over the subalgebra
of $\bool^c_0$ of those elements based on $\{x_\eta:\,\eta\in{}^{\omega>}
\lambda,\ \eta\neq \nu\}$.
\item[(2)] If $\eta$ is a  
branch of $f^\alpha$ hence necessarily $\eta\neq
\eta_\beta$ for $\beta\in Y'\cap \alpha$, $\xi<\blackz(\alpha)$, and
$w\subseteq \alpha\cap Y'$, is finite \then\ 
there is $k<\omega$ such that $\{\rho:\eta\restriction k
\vartriangleleft\rho\in {}^{\omega{>}}\lambda\}$ is disjoint to
$$
({}^{\omega>}\xi)\cup\bigcup\{N^\beta\cap {}^{\omega{>}}\lambda:\beta\in w,\
\beta +2^{\aleph_0}\leq\alpha\}\cup\bigcup\{{\rm supp}(s^\beta_n):\ n<\omega,
\beta\in w\}.
$$
\end{description}
\end{Fact}

\par \noindent
Proof. (1) Should be clear.

(2) Remember clauses (a),(c),(e) of Theorem \ref{1.7A}(2) and 
clause (1) of $\circledast_\alpha$ from stage A. 
\qquad \hfill$\qed_{\ref{2.4}}$

\relax From \ref{3.2} we can derive:
\medskip

\begin{Fact}
\label{3.4}
If $\xi<\blackz(\beta)$, $\beta<\alpha$, and $J\subseteq
{}^{\omega{>}}\lambda$ is finite, \then\  every element of
$\bool_\alpha$ which is based on $J\cup{}^{\omega>}\xi$ belongs
to $\bool_\beta$.
\end{Fact}

\par \noindent
Proof. We now prove by induction on
$\gamma\in[\beta,\alpha]$ that $[x\in \bool_\gamma\setminus \bool_\beta
\Rightarrow {\rm supp} (\lambda) \setminus {}^{\omega>}\xi$ is infinite].
For $\gamma=\beta$ this is empty, and for $\gamma$ limit it follows as
$\bool_\gamma=\cup\{\bool_\xi:\xi<\gamma\}$. For $\gamma+1\leq \alpha$,
let $x$ be a counterexample;\wolog\ $x\notin \bool_\gamma$; if 
$\bool_{\gamma+1}=\bool_\gamma$ this is impossible so $a_\gamma$, 
$\langle (d^\gamma_n, s^\gamma_n):n<\omega\rangle$ 
are well defined. Now $x$ is 
necessarily of the form $y_0\cup (y_1\cap a_\gamma) \cap (y_2-a_\gamma)$
 where $y_0, y_1,y_2$ are disjoint members of $\bool_\gamma$. 
Clearly $y_1\cap a_\gamma\notin \bool_\gamma$ or $y_2- a_\gamma
\notin \bool_\gamma$ so \wolog\ the former (otherwise use $-x$ which 
also $\in \bool_{\gamma+1} \setminus \bool_\gamma$ and has the same 
support). We can (by \ref{3.2}) find $n$ \st\ 
$J^*=\{\rho: \eta_\gamma\rest n 
\triangleleft \rho \in {}^{\omega>}\lambda\}$ 
is disjoint to $J$ and to $\supp (y_1)$. As $a_\gamma \cap 
\bigcup\limits_{\ell\leq n} d^\gamma_\ell \in \bool_\gamma$ clearly 
$y_2\cap a_\gamma-\bigcup\limits_{\ell\leq n} d^\gamma_\ell \notin 
\bool_\gamma$. As $\langle d^\gamma_m:m<\omega\rangle$ is a maximal
antichain of $\bool^c_0$, we can find $m\in (n,\omega)$ \st\ 
$d_m^\gamma\cap (y_2\cap a_\gamma-
\bigcup\limits_{\ell\leq n} d^\gamma_\ell)>0_{\bool_{\gamma+1}}$ hence 
$y_2\cap d^\gamma_m>0_{\bool_{\gamma+1}}$ and $s^\gamma_n$ has support 
$\subseteq J^*$ whereas $y_1,d^\gamma_m,x$ has 
support disjoint to $J^*$. But $y_1\cap d^\gamma_m\cap x=
(y_1\cap a_\gamma)\cap d^\gamma_m=(y_1\cap d^\gamma_m)\cap s^\gamma$ easy
contradiction.

\hfill $\qed_{\ref{3.4}}$
\medskip

\begin{Notation}
\label{3.5}
{\rm
\begin{description}
\item[(1)] Let $\bool^\xi$ be the set of $a\in \bool^c_0$
based on ${}^{\omega>}\xi$.
\item[(2)] For $x\in \bool^c_0$, $\xi<\lambda$ let $\rpr_\xi(x)=
\cap\{a\in \bool^\xi:x\leq a\}$.
\item[(3)] For $\xi\leq \lambda$ and $\nu\in {}^{\omega{>}}\xi$ let
$\bool^\xi_\nu$ be the set of $a\in \bool^c_0$ based on $J_{\xi,\nu}=:\{\rho\in
{}^{\omega{>}}\xi: \neg(\nu\vartriangleleft\rho)\}$. For $x\in \bool^c_0$ let
$\rpr_{\xi,\nu}(x)=\bigcap\{a\in \bool^\xi_\nu:x\leq a\}$.
\item[(4)] For $\gamma<\alpha^*$ let $\bool_{\langle\gamma\rangle}=\langle\{
x_\eta:\ \eta\in{}^{\omega>}\blackz(\gamma)\}\cup\{a_\beta:\beta\in\gamma
\cap Y'\} \rangle_{\bool_0^c}$.
\item[(5)] For $I\subseteq {}^{\omega >}\lambda$ and
$w\subseteq \alpha^\ast$ let

$\bool(I,w)=\langle\{x_\eta:\ \eta\in I\}\cup
\{a_\beta:\beta\in w\cap Y'\} \rangle_{\bool_0^c}$ and for $x\in\bool^c_0
,\xi\leq\lambda$
we let
$$
\rpr_{\xi,w}(x)=\cap\{y\in\langle \bool^\xi\cup\{x_\nu:\ \nu\in w\}
\rangle_{\bool^c_0}: x\leq y\}
$$
\item[(6)] For $\xi<\lambda$ let $\bool_{[\xi]}=\langle\{x_\eta:\ \eta\in
{}^{\omega>}\xi\}\cup\{a_\beta:\blackz(\beta)\leq\xi$ and
$\beta\in Y'\}\rangle_{\bool^c_0}$.
\item[(7)] For $J\subseteq {}^{\omega>}\lambda$ and $\xi\leq\lambda$ let 
$\rpr_{\xi,J} (\lambda)=\cap\{y\in\bool^\xi_J:x\leq y\}$ where 
$\bool^\xi_J=\langle\bool^\xi\cup\{x_\nu:\nu\in J\rangle_{\bool^c_0}$, 
when well defined.
\end{description}
}
\end{Notation}
\medskip

\begin{Fact}
\label{3.6}
\begin{description}
\item[(1)] For $\xi<\lambda$, $\bool^\xi$ is a complete Boolean subalgebra of
$\bool^c_0$. For $\xi<\lambda$ and $\nu\in {}^{\omega{>}}\xi$,
$\bool^\xi_\nu$ is a
complete subalgebra of $\bool^c_0$.
\item[(2)] If $\xi<\lambda$ and $x\in \bool^c_0$ \then\ 
$\rpr_\xi(x)$ is well defined
and belongs to $\bool^\xi$. Similarly, if $\xi<\lambda$, $\nu\in
{}^{\omega{>}}\xi$ and $x\in \bool^c_0$ \then\ 
$\rpr_{\xi,\nu}(x)$ is well
defined and belongs to $\bool^\xi_\nu$.
\item[(3)] If $\xi_0\leq\xi_1<\lambda$, $x\in \bool^c_0$ then
$\rpr_{\xi_0}(\rpr_{\xi_1} (x))=\rpr_{\xi_0}(x)$.
\item[(4)] If $\xi<\lambda$ and $w\subseteq \alpha^*$ is finite
\then\  the function $x \mapsto \rpr_{\xi,w}(x)$
is well defined for $x\in \bool^c_0$ and the value is in $\langle \bool^\xi\cup
\{a_\alpha:\alpha\in w\}\rangle_{\bool^c_0}$, of course which is a complete
subalgebra of $\bool^c_0$.
\item[(5)] If $\xi<\lambda$ and $\nu\in {}^{\omega{>}}\xi$ and $x\in \bool^c_0$
\then\  $\rpr_{\xi,\nu}(\rpr_\xi(x))=\rpr_{\xi,\nu}(x)$.
If in addition $\xi_0<\xi$ and
$\nu\notin {}^{\omega{>}}(\xi_0)$ \then\ 
$\rpr_{\xi_0}(x)=\rpr_{\xi_0}(\rpr_{\xi,\nu}(x))$.
\item[(6)] $\bool_{[\xi]}\subseteq \bool^\xi$ and if $\xi<\blackz
(\alpha)$ \then\  $\bool_{[\xi]} \subseteq \bool_\alpha$ 
and $\bool(I,w)\subseteq \bool_{\alpha^*}$
\end{description}
\end{Fact}

\par \noindent
Proof. Easy.

\begin{Fact}
\label{3.7}
\begin{description}
\item[(1)] For $x\in \bool_{\alpha^*}$, $\xi<\lambda$, the element 
$\rpr_\xi(x)$
belongs to $\bool_{[\xi]}$.
\item[(2)] If $x\in \bool_{\alpha^*}$, $\xi<\lambda$ and $J\subseteq {}^{\omega>}
(\xi+1)$ not necessarily finite, \then\  the element
$\rpr_{\xi,J}(x)$ belongs to $\langle \bool_{[\xi]}\cup \{x_\nu: \nu\in
J\}\rangle_{\bool^c_0}$.
\item[(3)] Like part (2) but $J\subseteq {}^{\omega>}\lambda$ (and not necessarily 
$J\subseteq {}^{\omega>}(\xi+1)$) and $J$ is finite 
\end{description}
\end{Fact}

\par \noindent
Proof. (1) We prove this for $x\in \bool_\alpha$, by induction on $\alpha$ (for
all $\xi$). Note that
$$
\boxdot\quad \rpr_\xi( \bigcup_{\ell<n} x_\ell)=\bigcup_{\ell<n} \rpr_\xi(x_\ell) 
\quad \mbox{ for } \quad x_0,\ldots, x_{n-1} \in \bool^c_0.
$$
\medskip

\noindent \underline{Case i:}\quad $\alpha=0$, or even just
$(\forall\beta<\alpha)[\blackz(\beta)\leq\xi]$.

\noindent Easy. Clearly we can find $\sigma$, $y_\ell$, $\nu_k$
($\ell< n$, $k<m$) such that
$x=\sigma(y_0,\ldots,y_{n-1}$, $x_{\nu_0},\dots,x_{\nu_{m-1}})$, 
where
$\sigma$ is a Boolean term, $y_\ell\in \bool_{[\xi]}$, $\nu_\ell\in
{}^\omega\lambda\setminus {}^{\omega>}\xi$; by the remark above 
\wolog\ 
$x=\bigcap\limits_{\ell<n+m} s_\ell$, where $s_\ell\in 
\{y_\ell,1-y_\ell\}$
when 
$\ell<n$, and $s_\ell\in\{x_{\nu_{\ell-n}},1-x_{\nu_{\ell-n}}\}$ 
when $n
\leq\ell<n+m$, and the sequence $\langle x_{\nu_0},\ldots,x_{\nu_{n-1}}
\rangle$ is without repetitions. Now by \ref{3.3} clearly 
$\rpr_\xi(x)=\bigcap\limits_{
\ell<n} s_\ell$ which belongs to $\bool_{[\xi]}$;
\medskip

\noindent \underline{Case ii:}\quad $\alpha$ is limit.

\noindent Trivial as $\bool_\alpha=
\bigcup\limits_{\beta<\alpha}\bool_\beta$.
\medskip

\noindent \underline{Case iii:} $\alpha=\beta+1$.

\noindent By the induction hypothesis \wolog\  
$x\not\in \bool_\beta$ hence $\beta\in Y'$. As $x\in
\bool_\alpha$ there are disjoint $\ebool_0$, $\ebool_1$, 
$\ebool_2\in \bool_\beta$
such that $x=\ebool_0
\cup(\ebool_1\cap a_\beta)\cup (\ebool_2-a_\beta)$. 
It suffices to prove that $\rpr_\xi(
\ebool_0)$, $\rpr_\xi(\ebool_1\cap a_\beta), 
\rpr_\xi(\ebool_2-a_\beta)\in \bool_{[\xi]}$; the
first holds by the induction hypothesis and \wolog\ 
~we concentrate on
the second. Remembering
clause $\circledast$(1) of stage (A), by \ref{3.2} 
applied to $\bool_\alpha$, $\ebool_1$
we have: there are
$\xi_0<\blackz(\beta)$ and $k<\omega$ such that $\ebool_1$ is based on $J
\stackrel{\rm def}{=}{}^{\omega>}\lambda\setminus\{\rho:\eta_\beta
\restriction k\vartriangleleft\rho\in{}^{\omega>}\lambda\}$. Now
\wolog\  each
$d^\beta_n$ ($n<\omega$) is based on ${}^{\omega>}\xi_0$ 
(recall clause  $\circledast_\beta$ (2) of
stage A) and ${}^{\omega>}
\xi_0\subseteq J$ (this holds if $\eta_\beta\restriction k\notin
{}^{\omega{>}}\xi_0$, and as $\eta_\beta$ is increasing with limit
$\blackz(\beta)$ this is easy to obtain). By Case i, we can assume that
$\xi<\blackz(\beta)$ hence (as we can increase $k$ and $\xi_0$)
\wolog\ ~$\xi<\xi_0$, and
by the induction hypothesis and \ref{3.6}(3),(5), 
letting $\nu=:\eta_\beta
\restriction k$, it suffices to prove $\rpr_{\xi_0,\nu}
(\ebool_1\cap a_\beta)\in
\bool_{\beta}$. As $m<\omega\ \Rightarrow\ a_\beta\cap d^\beta_m\in
\bool_{[\blackz(\beta)]}$ and $\boxdot$ above, 
\wolog\ ~$\ebool_1\cap d^\beta_m=0$ for $m<k$.
Now clearly for proving $\rpr_{\xi_0,
\nu}(\ebool_1\cap a_\beta)=\ebool_1$ 
it is enough to show, for each $m<\omega$, that
$\rpr_{\xi} (\ebool_1\cap d^\beta_m\cap \bools^\alpha_m)=
\ebool_1\cap d^\beta_m$ as $\langle
d^\beta_n:n<\omega\rangle$ is a maximal antichain of $\bool^c_0$ and as
$a_\beta\cap d^\beta_m=\bools^\beta_m$ both by $\circledast_\beta (2)$.
If $m<k$ then
$\ebool_1\cap d^\beta_m=0$ so this is trivial. 
If $m\geq k$ this holds because
$d^\beta_m, \ebool_1$ are based on $J$, 
${}^{\omega>}\xi_0\subseteq J$ and
$\bools^\beta_m$ is based on ${}^{\omega>}\lambda\setminus J$ and is 
$\bools^\beta_m >0$.

\par \noindent 
(2),(3) Same proof. \qquad\hfill $\qed_{\ref{3.7}}$
\smallskip

\begin{Lemma}
\label{3.8}
1) Suppose that $I$, $w$ satisfy:
\begin{description}
\item[$(*)_{I,w}$] $I\subseteq {}^{\omega>}\lambda$,
$w\subseteq\alpha^*\cap Y'$,
$I$ is closed under initial segments,
\[\alpha\in w\ \&\ n<\omega\quad\Rightarrow\quad\eta_\alpha\restriction n\in
I,\]
and for every $\alpha<\alpha^*$, if $\bigwedge\limits_{m<\omega}(
\eta_\alpha\restriction m\in I)$ then $s^\alpha_m, d^\alpha_m$ are based on
$I$ and belong to $\bool(I,w)$ and $\alpha\in w$; 
see Definition \ref{3.5}(5).
\end{description}
\Then\ for any countable 
$C\subseteq \bool_{\alpha^*}$ there is
a projection $h$ from $\langle \bool(I,w),C\rangle_{\bool_0^c}$
onto $\bool(I,w)$.
\par \noindent
2) If $(*)_{I,w}$ holds \then\ every member of $\bool(I,w)$ is based on $I$.
\par \noindent
3) We can add
\begin{description}
\item[(a)] if $a_\alpha\in C \setminus \bool(I,w)$,
and $\{d^\alpha_n: n< \omega\}\subseteq
\bool(I, w)$ and $\{\eta_\alpha\restriction n:
n< \omega\}\subseteq I$
\then\  $h(a_\alpha)$ has support $\subseteq{}^{\omega>}\zeta$ for
some $\zeta< \blackz(\alpha)$.
\item[(b)] if $\nu\in {}^{\omega>}\lambda$ and $x_\nu\in C$ 
\then\ 
$h(x_\nu) \in \{0, x_\nu, 1\}$
\item[(c)] if $c=\sigma(a_{\alpha_0}, \ldots, a_{\alpha_{k-1}}, b_0,
\ldots, b_{n-1})$ where $\sigma$ is a Boolean term
$\alpha_0,\ldots \alpha_{k-1} \in w$ and
$\blackz(\alpha_0)\leq\ldots
\leq\blackz(\alpha_{k-1})\leq\xi$ and $b_\ell\in \bool
(I, w)$ and ${\rm supp}
(b_\ell)\subseteq {}^{\omega>} \zeta$ for some 
$\zeta<\xi$ and
$a_{\alpha_\ell}\in C$ \then\  $h(c) = \sigma(a'_0,
\ldots, a'_{k-1}, b_0, \ldots, b_{n-1})$ (where $\ell< k\
\Rightarrow\ a'_\ell = h(a_{\alpha_\ell}))$ is in
$\bigcup\limits_{\varepsilon <\xi} \bool_{[\varepsilon]}$.
\end{description}
\end{Lemma}

\noindent\underline{Remark:}\quad In $(*)_{I,w}$ the last phrase can be
weakened to ``for some $m_\alpha<\omega$, for every $m\in [m_\alpha,\omega)$
the elements $s^\alpha_n, d^\alpha_m$ are based on $I$ (and belong to $\bool(I,w)$
and $(\alpha\in w)$''.
\medskip

\par \noindent
Proof. 1)  We can easily find $I(*)$, $w(*)$ such that
$C\subseteq \bool(I(\ast),
w(*))$, $w\subseteq w(*)\subseteq\alpha^*$, $|w(*)\setminus w|\leq\aleph_0$,
$I\subseteq I(*)\subseteq {}^{\omega>}
\lambda$, $I(*)$ is closed under
initial segments, $|I(*)\setminus I|\leq \aleph_0$,
and if $\alpha\in w(*)
\setminus w$, then $\bools^\alpha_m, d^\alpha_m\in \bool(I(*),w(*))$ hence 
$\{\eta_\alpha\rest m:m<w\}\subseteq I(*)$. Let
$w(*)\setminus w=\{\alpha_\ell:\ell<\omega\}$ for notational simplicity,
and we choose by induction on
$\ell$ a natural number $k_\ell<\omega$, such that the sets
\[\{\nu\in {}^{\omega>}\lambda:\nu\mbox{ appears in 
$\bools^{\alpha_\ell}_m$
for some }m\geq k_\ell\}\]
are pairwise disjoint and disjoint to $I$ (possible by the demand
$s^\alpha_m\in \langle x_\rho:\eta_\alpha\restriction
m\vartriangleleft \rho, \rho\in {}^{\omega>}\lambda\rangle_{\bool^c_0}$ in
clause $\circledast_\alpha$(2) of stage A in the beginning of the proof
of \ref{3.1} and $(*)_{I,w}$). First assume that 

$\boxdot$ \,for every $\alpha, \langle \supp 
(s^\alpha_m):m<\omega\rangle$ 
is a \sq\ of pairwise disjoint sets. 

\par \noindent
Now we can extend the identity map on $\bool(I,w)$ to a projection $h_0$ from
$\bool(I(*),w)$ onto $\bool(I,w)$ such that
\begin{description}
\item[(a)] if $\nu\in I(\ast)\setminus
I$ then $h_0(x_\nu) \in \{0, 1\}$
\item[(b)] if $\ell<\omega$, $m>k_\ell$, then 
$h_0(s^{\alpha_\ell}_m)=0$.
\end{description}
This is possible as $\bool(I(\ast), w)$ is generated by
 $\bool(I,w)\cup\{x_\nu: \nu\in I(*)\setminus I\}$
freely except the equations which  hold in $\bool(I, w)$ and 
$\boxdot$ above as $\langle d^{\alpha_\ell}_m:m\in (k_\ell,\omega)\rangle$
is a \sq\ of  pairwise disjoint elements.
Now we can define by induction on $\alpha\in (w(*)
\setminus w)\cup \{\alpha^*\}$ a projection $h_\alpha$ from
$\bool(I(*),w\cup (w(*)\cap\alpha))$ onto
$\bool(I,w)$ extending $h_\beta$ for any $\beta<\alpha$
satisfying $\beta\in (w(*)\setminus w)\cup\{0\}$. For $\alpha=0$ we have
defined it, for $\alpha=\alpha^*$ we get the desired
conclusion, and in limit stages
take the union. In successive stages there is no problem by the choice of
$h_0$, and of the $k_\ell$'s (and $\circledast$(2) of stage A).

If $\boxdot$ fails, we just define $h_\xi=h\rest
(\bool (I(*),w(*))\cap\bool^\xi)$ by induction on $\xi\leq \lambda$ such 
that (it is the identity on $\bool (I,w)\cap \rDom 
(h_\xi))$ and 
\begin{enumerate}
\item[$(a)'$] if $\nu\in I(*)\setminus I$ and $x_\nu \in \rDom
(h_\xi)$ then $h_\xi(x_\nu)\in \{0,1\}$
\item[$(b)'$] if $\ell<\omega,m>k_\ell$ and $h_\xi (d^{\alpha_\ell}_m)$
 is well defined then $s^{\alpha_\ell}_m \cap d^{\alpha_\ell}_m$ 
does not belong to the filter on $\bool^c_0$ generated by 
$\{d\in \rDom (h_\xi):h_\xi (d)=1\}$. 
\end{enumerate}
2) The proof of part 1) gives this. 
\par \noindent
3) Note that by clause $(a)'$, in clause $(b)'$, if 
$h_\xi (s^{\alpha_\ell}_m\cap d^{\alpha_\ell}_m)$ is well defined then it is 
$0_{\bool_0}$
This is possible by the choice of $\langle k_\ell:\ell<\omega\rangle$ and as 
$\langle d^{\alpha_\ell}_m:m<\omega\rangle$ is a \sq\ of pairwise 
disjoint elements of $\bool^c_0$.   
\qquad \hfill$\qed_{\ref{3.8}}$
\medskip

\begin{Claim}
\label{3.9}
If $\bool'$ is an uncountable subalgebra of 
$\bool_{\alpha^*}$ \then\ 
 there is an antichain $\{d_n:n<\omega\}
\subseteq \bool'$ such that for no
 $x\in \bool_{\alpha^*}$ do we have $x\cap d_{2n}=0$, $x\cap
 d_{2n+1}=d_{2n+1}$ for every $n$, provided that
\begin{description}
\item[(*)] no single countable $I\subseteq {}^{\omega >}\lambda$
is a support for every $a\in \bool'$.
\end{description}
\end{Claim}

Proof:  We choose by induction on $\alpha<\omega_1$, $d_\alpha$,
$I_\alpha$, such that:
\begin{description}
\item[(i)]   $I_\alpha\subseteq {}^{\omega>}\lambda$ is countable,
closed under initial segments
\item[(ii)]  $\bigcup\limits_{\beta<\alpha}
I_\beta\subseteq I_\alpha$ and
for $\alpha$ limit, equality holds,
\item[(iii)] $d_\alpha\in \bool'$ is based on $I_{\alpha+1}$
but not on $I_\alpha$.
\end{description}
There is no problem doing this as we are assuming (*).

By clause (iii), for each $\alpha$ there are a non zero
$s^0_\alpha\in\langle x_\eta:
\eta\in I_\alpha\rangle_{\bool^c_0}$
and non-zero $s^1_\alpha$, $s^2_\alpha\in
\langle x_\eta:\eta\in I_{\alpha+1}\setminus
I_\alpha\rangle_{\bool^c_0}$ such
that $s^1_\alpha\cap s^2_\alpha=0$, 
$s^0_\alpha\cap s^1_\alpha\leq
d_\alpha$, $s^0_\alpha\cap s^2_\alpha\leq 1-d_\alpha$.

By Fodor's lemma, as we can replace $I_\alpha$ by
$I_{h(\alpha)}$ if
$h:\omega_1\ \rightarrow\ \omega_1$ is increasing continuous,
\wolog, $s^0_\alpha=s^0$ (i.e., does not depend on
$\alpha$). For each $\alpha$ there is $n(\alpha)<\omega$ such that
$$
s^0=s^0_\alpha\,\in\, \langle x_\eta:\eta\in I_\alpha\cap {}^{n(\alpha)\geq}
\lambda\rangle_{\bool^c_0},
$$
$$
s^1_\alpha,s^2_\alpha\, \in\, \langle x_\eta:\eta\in (I_{\alpha+1}\setminus
I_\alpha)\cap {}^{n(\alpha)\geq}\lambda\rangle_{\bool^c_0}.
$$
Again, by renaming \wolog\ ~$n(\alpha)=n(*)$ for every $\alpha$. For $n<
\omega$ let $d^n=d_n-\bigcup\limits_{\ell<n} d_\ell$, $s^n=
s^0\cap
\bigcap\limits_{\ell<n} s^2_\ell\cap s^1_n$, 
so easily $d^n\in \bool'$, $\langle
d^n:n<\omega\rangle$ is an antichain, $s^n\leq d^n$ and 
$s^n\in\langle
x_\eta:\ \eta\in {}^{n(*)\geq}\lambda\rangle_{\bool^c_0}$
and by the choice of $\bool_0$
easily $0<s^n$.
 Suppose $x\in \bool_\alpha$ satisfies:
for each $n<\omega$, we have $x\cap d^{2n}=0$, $x\cap d^{2n+1}=d^{2n+1}$. Then
for $n<\omega$, $x\cap s^{2n}=0$, $x\cap 
s^{2n+1}=s^{2n+1}$. But by
\ref{3.8}(1) (for $I={}^{n(*)\geq}\lambda$,
$w=\emptyset$ and $C=\{x\}$), there
is such $x$ in $\langle x_\eta:\eta\in{}^{n(*)\geq}\lambda\rangle_{\bool^c_0}$,
an easy contradiction. \qquad \hfill$\qed_{\ref{3.9}}$
\medskip

Hence we have proved in particular that for every $\aleph_1$-
compact $\bool'\subseteq \bool_{\alpha^*}$, some countable
$I\subseteq {}^{\omega >}\lambda$ supports
every $x\in \bool'$.

\begin{Claim}[Crucial Claim]
\label{3.10}
No infinite subalgebra $\bool'$ of $\bool_{\alpha^*}$ is $\aleph_1$-compact.
\end{Claim}
 \par \noindent
Proof. Suppose that there is such subalgebra, and let $\xi$ be minimal such
that there is an infinite $\aleph_1$--compact $\bool'\subseteq
\bool_{[\xi]}$. The proof is broken into five parts.
\medskip

\noindent \underline{Part I}\quad If
\begin{description}
\item[(a)] $\bool'\subseteq \bool_{\alpha^*}$ is $\aleph_1$-compact and
infinite (subalgebra)
\item[(b)] $\bool'\subseteq \bool_{[\xi]}$,
\end{description}
then
\begin{description}
\item[(c)] for every $\zeta<\xi$, finite $J\subseteq {}^{\omega>}\lambda$
and $x\in
\bool'\setminus\{y:\{z\in \bool':z\leq
y\}$ is finite$\}$, there is $x_1\in \bool'$, $x_1\leq x$ such that for no
$y\in \langle\bool_{[\zeta]}\cup\{x_\nu:\nu\in J\}\rangle_{\bool^c_0}$, do we have $y\cap x=x_1$.
\end{description}
So toward contradiction assume $\bool'$ satisfies (a) and (b), but it
fails (c) for $\zeta<\xi$, a finite $J\subseteq {}^{\omega>}\lambda$ 
and $x\in \bool'$, hence $\{y:\ y\leq x,\ y\in
\bool'\}$ is infinite. 
So for every $z\in \bool'$, there is $g(z)\in
\langle \bool_{[\zeta]}\cup \{x_\nu:\nu\in J\}\rangle_{\bool^c_0}$.  
such that $g(z)\cap x=z\cap x$ (otherwise we can use
$x_1=z\cap x$). Let $\bool^a$ be the Boolean subalgebra of 
$\langle \bool_{[\zeta]}\cup \{x_\nu:\nu\in J\}\rangle_{\bool^c_0}$ 
generated by $\{g(z):z\in \bool'\}$, so $z\in \bool'\ \Rightarrow\
z\cap x\in \bool^a$. Clearly
\[\{y\in \bool':y\leq x\}=\{t\cap x:t\in \bool^a\}.\]
Let $x^*=\rpr_{\zeta,J}(x)$ (it is in $\bool_{[\zeta]}$ 
by \ref{3.7}(1) if $J=\emptyset$, \ref{3.7}(3) otherwise), and let
\[\bool^b=\{t\cap x^*:\ t\in \bool^a\}\cup\{t\cup (1-x^*):t\in \bool^a\}.\]
Clearly $\bool^b$ is a subalgebra of 
$\langle \bool_{[\zeta]}\cup \{x_\nu:\nu\in J\}\rangle_{\bool^c_0}$, 
and $1-x^*$
is an atom of
$\bool^b$. Now $\bool^b$ is infinite, why?  there are distinct $x_n\leq x$ in
$\bool'$ (for $n<\omega$), so $g(x_n)\in \bool^a$ and hence
$g(x_n)\cap x^*\in \bool^b$. As $x\leq x^*$ and
\[n\neq m\quad \Rightarrow\quad g(x_n)\cap x=x_n\cap x=x_n\neq x_m=x_m\cap
x=g(x_m)\cap x,\]
clearly $[n\neq m\ \Rightarrow\ g(x_n)\cap x^*\neq g(x_n)\cap x^*]$
so $\bool^b$ is really infinite. We
shall prove that $\bool^b$ is $\aleph_1$-compact, thus contradicting the choice
of $\xi$. Let $d_n\in \bool^b$ be pairwise disjoint, and we would like 
 to find $t\in
\bool^b$ satisfying $t\cap d_{2n}=0$, $t \cap d_{2n+1}= d_{2n+1}$ for $n<\omega$.
Clearly \wolog\ ~$d_n \leq x^*$ (as $1-x^*$ is an atom of $\bool^b$). So $d_n=
t_n\cap x^*$ for some $t_n\in \bool^a$, hence easily $t_n \cap x\in \bool'$
so for some $x_n\in \bool'$, $x_n\leq x$ and $t_n\cap x=x_n\cap x=x_n$.
So $x_n=g(x_n) \cap x$.

For $n\neq m$,
\[x_n\cap x_m=(t_n\cap x)\cap (t_m\cap x)\leq (t_n\cap x^*)\cap (t_m\cap
x^*)= d_n\cap d_m=0.\]
As $\bool'$ is $\aleph_1$-compact there is $y\in \bool'$ satisfying
$y\cap x_{2n}
=0$, $y\cap x_{2n+1}=x_{2n+1}$ for $n<\omega$. Now $g(y)$, $d_n$, $t_n$
belong to $\langle \bool_{[\zeta]}\cup 
\{x_\nu:\nu\in J\}\rangle_{\bool^c_0}$ and 
(as $x_n\leq x \leq x^*$ and $d_n=t_n\cap x^*$,
$t_n\cap x=x_n$):
\begin{description}
\item[(i)] $g(y)\cap d_{2n}\cap x=g(y)\cap t_{2n}\cap x=g(y)\cap x_{2n}\cap x
=y\cap x_{2n} \cap x=0$,
\item[(ii)] $g(y)\cap d_{2n+1}\cap x=g(y)\cap t_{2n+1}\cap x=g(y)\cap
x_{2n+1}\cap x= y\cap x_{2n+1}\cap x=x_{2n+1}\cap x=t_{2n+1}\cap x=d_{2n+1}
\cap x$.
\end{description}
Now, by the definition of $x^*=\rpr_{\zeta,J}(x)$,
\[s\in \langle \bool_{[\zeta]}\cup 
\{x_\nu:\nu\in J\}\rangle_{\bool^c_0} 
\& s\cap x=0\quad\Rightarrow\quad s\cap x^*=0\]
(as $1-s\in \langle\langle \bool'_{[\zeta]}\cup 
\{x_\nu:\nu\in J\}\rangle_{\bool^c_0}$ and by the left side 
$x\leq 1-s$), hence by clause (i) (for 
$s=g(y)\cap d_{2n}$):
\begin{description}
\item[(iii)] $g(y)\cap d_{2n}\cap x^*=0$.
\end{description}
Also, by the definition of $x^*=\rpr_{\zeta,J}(x)$,
\[s_1,s_2\in \langle 
\bool_{[\zeta]}\cup \{x_\nu:\nu\in J\}_{\bool^c_0}
\ \& \ s_1\cap x=s_2\cap
x\quad\Rightarrow\quad s_1 \cap x^*=s_2\cap x^*\]
(as $s_1-s_2\in \bool_{[\zeta]}\cup 
\{x_\nu:\nu\in J\}\rangle_{\bool^c_0}$
and by the left side $x\leq 1-
(s_1-s_2)$ hence as above  $x^*\leq 1-(s_1-s_2)$ 
and similarly $x^*
\leq 1-(s_2- s_1)$). Hence by clause (ii)
\begin{description}
\item[(iv)] $g(y)\cap d_{2n+1}\cap x^*=d_{2n+1}\cap x^*$.
\end{description}
But $d_n\leq x^*$, so from (iii) and (iv), $(g(y)\cap x^*)\cap d_{2n}=0$,
$(g(y)\cap x^*)\cap d_{2n+1}=d_{2n+1}$, and $g(y)\in \bool^a$, hence $g(y)\cap
x^*\in \bool^b$. So $\bool^b$ is $\aleph_1$--compact and this contradicts the
minimality of $\xi$, hence we finish proving Part I.
\medskip

\noindent\underline{Part II}:\quad If $\bool^1\subseteq 
\bool_{\alpha^*}$ is $\aleph_1$--
compact, $\bool^1 \subseteq \bool^2$, $\bool^2=\langle
\bool^1\cup\{z\}\rangle_{\bool^2}$ then $\bool^2$ is $\aleph_1$--compact.

\noindent The proof is straightforward. [If $d_n\in \bool^2$ are pairwise
disjoint, let $d_n=d^0_n \cup(d^1_n\cap z)\cup (d^2_n-z)$ for some 
disjoint $d^0_n,d^1_n$, $d^2_n\in
\bool^1$. As $\bool_{\alpha^*}$ satisfies the c.c.c. also $\bool^1$ 
satisfies the c.c.c., hence being $\aleph_1$-compact, is complete. Now the 
each of the sets each $\idealJ_\ell(\ell<3)$ is an ideal of $\bool^1$
and their union 
$\idealJ_0\cup \idealJ_1\cup \idealJ_2$ is a dense subset of 
$\bool^1$ where $\idealJ_\ell= \{
x\in \bool^1:x>0$ satisfies $\ell=0\Rightarrow\bool^2\models x\cap z=0$ and
$\ell=1\Rightarrow \bool^2\models x\leq z$ and $\ell=2\Rightarrow \bool^2
\models (\forall y) (0<y\leq x \ \& \ y\in \bool^1\Rightarrow y\cap z\neq
0\neq y-z)\}$. As $\bool^1$ is complete \wolog\  
$d^\ell_m\in \idealJ_0 \cup \idealJ_1\cup 
\idealJ_2$ for $m<\omega,\ell<3$. Also there is a 
maximal antichain  $\langle x_n:n<\gamma\leq \omega\rangle_{\bool^1}$ 
of $\bool^1$ consisting of 
elements of this family. Similarly \wolog\ for each 
$n$ we have  $x_n\leq d^\ell_m$ for 
some $m<\omega_1$ and $\ell<3$; or $x \cap d^\ell_m=0$
 for every $m$.  \Wlog\ 
$d^1_n\neq 0 \Rightarrow d^1_n\in \idealJ_2$ and 
$d^2_n\neq0\Rightarrow d^2_n\in \idealJ_2$ and necessarily 
$d^0_n \cap (d^1_m\cup d^2_m)=0$ for $n,m<\omega$. 
Now, necessarily $d^0_n \cap d^0_m=0$ for $n\neq m$ and  \wolog, 
$d^1_n\cap d^1_m=0$ for $n\neq m$  --- otherwise
replace them by $d^1_n-\bigcup\limits_{\ell<n} d^1_\ell$; 
similarly $d^2_n
\cap d^2_m=0$, for $n\neq m$. So, for $\ell=0,1,2$, there is $y^\ell\in \bool^1$
such that for every $n<\omega$ we have:
\[y^\ell\cap d^\ell_{2n}=0,\quad y^\ell\cap d^\ell_{2n+1}=d^\ell_{2n+1}.\]
Hence $y^0\cup (y^1\cap z-y')\cup (y^2\cap z-y')$ is a solution.]
\medskip

\noindent\underline{Part III}:\quad $\xi$ cannot be a successor ordinal.

\noindent Proof:  Let $\bool'$ satisfy clauses (a), (b) (hence (c)) of Part I.

Suppose toward contradiction that $\xi=\zeta+1$, and by \ref{3.9} there is
a countable $I\subseteq {}^{\omega>}\xi$ which supports every $a\in
\bool'$; \wolog, $I$ is closed under initial segments and, under those
demands, $|I\setminus {}^{\omega>}\zeta|\leq\aleph_0$ is
minimal. Now, by applying Part
I we get
\begin{description}
\item[$\boxdot$] for every finite $J\subseteq {}^{\omega>}\lambda$, and $x\in \bool'$ for which
$\{y\in \bool':y\leq x\}$ is infinite, there is $x_1\in \bool'$, $x_1\leq x$
such that for no $y\in\langle \bool_{[\zeta]}\cup \{x_\eta:\eta\in
J\}\rangle_{\bool^c_0}$ do we have $y \cap x=x_1$.
\end{description}
Now, $I\setminus {}^{\omega>}\zeta$ is infinite.
[Why? Otherwise let $\bool''=
\langle \bool'\cup\{x_\eta:\eta \in I\setminus {}^{\omega>}\zeta\}
\rangle_{\bool^c_0}$; it is infinite and $\aleph_1$-
compact by Part II, and we shall
we apply Part I to it. Let $k=|I\setminus {}^{\omega{>}}\zeta|$
and let $I
\setminus {}^{\omega>}\zeta=\{\eta_0,\dots,\eta_{k-1}\}$
and for $u\subseteq \{0,\dots,k-1\}$, let
\[x_u\stackrel{\rm def}{=}\bigcap\{x_{\eta_\ell}:\ell\in u\}
\cap\bigcap\{1-
x_{\eta_\ell}:\ell<k,\quad \mbox{ and } \quad \ell\notin u\}.\]

So $x_u\in \bool''$, $1=\bigcup\{x_u:u\subseteq \{0,\dots,k-1\}\}$, hence for
some $u$, $\{y\in \bool'':y\leq x_u\}$ is infinite;
now $\zeta$, $x_u$ contradict
the conclusion of Part I.]

As $\bool'$ is $\aleph_1$--compact, for any $x\in \bool'$
such that $\{y\in \bool':y
\leq x\}$ is infinite, $x$ can be splitted in $\bool'$ to two elements
satisfying the same, i.e., $x=x^1\cup x^2$, $x^1\cap x^2=0$, $\{y\in \bool':y
\leq x^\ell\}$ is infinite for $\ell=1,2$. Let $I\setminus {}^{\omega>}
\zeta=\{\eta_\ell:\ell<\omega\}$, so we can find pairwise disjoint 
$\ebool_n\in
\bool'$ such that $\{y\in \bool':y\leq \ebool_n\}$ is infinite.
Now, by $\boxdot$ above, for
each $n$ we can find $d_{2n}$, $d_{2n+1}$ satisfying $\ebool_n=d_{2n}\cup
d_{2n+1}$, $d_{2n}\cap d_{2n+1}=0$ and such that for no $y\in\langle
\bool_{[\zeta]}\cup\{x_{\eta_\ell}:\ell<n\}\rangle_{\bool^c_0}$
do we have $y\cap(d_{2n}\cup d_{2n+1})=d_{2n+1}$.

Since $\bool'$ is $\aleph_1$--compact there is $y\in \bool'$ such that $y\cap
(d_{2n}\cup d_{2n+1})=d_{2n+1}$ for every $n<\omega$. As $y\in \bool'$ clearly
$y\in \bool_{[\xi]}=\bool_{[\zeta+1]}$, and $y$ is based on 
$\{x_\nu:\nu\in {}^{\omega>}\zeta\}\cup\{
x_{\eta_\ell}:\ell<\omega\}$, so by \ref{3.7}(2) we have  
$y'=\rpr_{\zeta,\{\eta_\ell:\ell<\omega\}} (x)$ belong to $\langle
\bool_{[\zeta]}\cup\{x_{\eta_\ell}:\ell<\omega\}\rangle_{\bool^c_0}$.
Hence $y'\in \langle \bool_{[\zeta]}\cup\{x_{\eta_\ell}:\ell<n\}
\rangle_{\bool^c_0}$ for some $n<
\omega$. This is a contradiction to $y'\cap(d_{2n}\cup d_{2n+1})=d_{2n+1}$ 
which holds as by the choice of $d_{2n},d_{2n+1}$, so $y\cap d_{2n}=0$, 
$y\cap d_{2n+1}=d_{2n+1}$ and $d_{2n},d_{2n+1}\in \langle\bool_{[\zeta]}\cup 
\{x_{\eta_\ell}:\ell<\omega\}\rangle_{\bool_0^c}$ so $y'\cap d_{2n}=0, 
y'\cap d_{2n+1}=d_{2n+1}$.
\medskip

\noindent\underline{Part IV}:\quad Let $\bool'$ satisfy
clauses (a), (b) of Part
I (and hence clause (c) too). By \ref{3.9}, for some countable
$I\subseteq
{}^{\omega>}\xi$, every $b\in \bool'$ is based on $I$.
By Part III, $\xi$ is not
a successor ordinal and trivially it is not zero hence $\xi$ is a
limit ordinal. Now by \ref{3.5}(6) (i.e. the definition of $\bool_{[\zeta]}$
for $\zeta\leq \lambda$) for no $\zeta< \xi$ is
$I\subseteq {}^{\omega >}\zeta$,
hence necessarily $\cf(\xi)=\aleph_0$. Let
\[Fi(\bool')=\{x\in \bool':
\{y\in \bool':y\leq x\}\mbox{ is finite}\}.\]  
Next we shall show:
\begin{description}
\item[$(**)$] for some finite $w^\ast\subseteq\{\alpha<\alpha^\ast:
\blackz(\alpha)=\xi\}$ and
$x^*\in \bool'\setminus Fi(\bool')$, for every $y\leq x^*$
 from $\bool'$, for some $z\in
\langle\cup_{\zeta<\xi}\bool_{[\zeta]}\cup\{a_\alpha:\alpha\in w^\ast\}
\rangle_{\bool^c_0}$ we have $z\cap x^*=y$.
\end{description}
Suppose $(**)$ fails and we choose by induction on $n<\omega$, $x_n$, $y_n$,
$w_n$ such that:
\begin{description}
\item[(i)]   $x_n\in \bool'$, and $m<n\ \Rightarrow\ x_m\cap x_n=0$,
\item[(ii)]  $1-\bigcup\limits_{i<n} x_i \not\in Fi(\bool')$,
\item[(iii)] $w_n\subseteq \{\alpha:\blackz(\alpha)=\xi\}$ is finite,
\item[(iv)]  $w_n\subseteq w_{n+1}$,
\item[(v)]   $y_n\leq x_n$ and $y_n\in \bool'$,
\item[(vi)]  for no $z\in\langle\bigcup\limits_{\zeta<\xi} \bool_{[\zeta]}\cup
\{a_\alpha:\alpha\in w_n\}\rangle_{\bool^c_0}$ do we have $z\cap x_n=y_n$.
\end{description}
For $n=0$ we have $1\notin Fi(\bool')$, hence (ii) is satisfied.

\noindent For each $n$ let $w_n$ be a finite subset of $\{\alpha:\blackz(
\alpha)=\xi\}$ extending $\bigcup\limits_{\ell<n} w_\ell$ such that for
every $\ell<n$, $x_\ell$, $y_\ell\in\langle\bigcup\limits_{\zeta<\xi}
\bool_{[\zeta]}\cup \{a_\alpha:\alpha\in w_n\}\rangle_{\bool^c_0}$, it exists
by \ref{3.5}(6). Then, as $1-\bigcup\limits_{\ell<n} x_\ell\notin
Fi(\bool')$, and as $\bool'$ is $\aleph_1$--compact, there is $x_n\leq
1-\bigcup\limits_{i<n} x_i$ satisfying $x_n\in
\bool'$ such that $1-\bigcup\limits_{\ell\leq n} x_\ell\notin Fi(\bool')$ and $x_n
\notin Fi(\bool')$. Now, as (**) fails, necessarily $w_n$, $x_n$ do not satisfy
the requirements on $w^*$, $x^*$ in $(**)$, so there is $y_n\in \bool'$,
$y_n\leq x_n$ such that for no $z\in\langle\bigcup\limits_{\zeta<\xi}
\bool_{[\zeta]}\cup
\{a_\alpha:\alpha\in w_n\}\rangle_{\bool^c_0}$ 
do we have $z\cap x_n=y_n$.
So we can carry the definition. As $\bool'$ is $\aleph_1$--
compact, for some $z^\ast\in \bool'$ we have
$z^\ast\cap x_n=y_n$ for every $n$.

As $z^\ast\in \bool'$ and $\bool'\subseteq \bool_{[\xi]}$,
for some finite $w^\ast
\subseteq\{\alpha<\alpha^\ast:\blackz(\alpha)=\xi\}$ we have
\[z^\ast=\sigma(\ldots,a_\alpha,\ldots,\ldots,b_\ell,\ldots)_{\alpha\in
w^*,\ell<n}\in\langle\bigcup\limits_{\varepsilon<\xi}
\bool_{[\varepsilon]}\cup\{a_\alpha:
\alpha\in w^\ast\}\rangle_{\bool^c_0}\]
where $\sigma$ is a Boolean term, and $\ell<n \ \Rightarrow\ b_\ell\in
\bigcup\limits_{\epsilon <\xi} \bool_{[\epsilon ]}$.
As $w^\ast$ is finite, for some $n(\ast)<\omega$ we have $w^\ast\cap
(\bigcup\limits_{n<\omega} w_n)\subseteq w_{n(\ast)}$.

Let $k^*<\omega$ be
such that there are no repetitions in 
$\langle\eta_\alpha\restriction k^*:
\alpha\in w_{n(*)+1}\rangle$ and $k^\ast>n(*)$.
Let $\zeta<\xi$ be such that: ${\rm supp}(
d^\alpha_n)\subseteq {}^{\omega>}\zeta$ for $\alpha\in w_{n(\ast)+1}\cup
w^\ast$, $n<\omega$ and 
${\rm supp}(\bools^\alpha_k)\subseteq{}^{\omega{>}}\zeta$
for $\alpha\in w_{n(*)+1}\cup w^*$, $k<k^*$, and
\[x_n,y_n\in\langle \bool_{[\zeta]}
\cup\{a_\alpha:\alpha\in w_{n(\ast)+1}\}\rangle_{\bool^c_0}\]
for $n<n(\ast)+1$ and $z^\ast\in\langle 
\bool_{[\zeta]}\cup\{a_\alpha:\alpha\in w^\ast\}\rangle_{\bool^c_0}$

 We shall now apply \ref{3.8}
with $I,w,C$ there standing for
\[I'=\{\eta:\eta\in{}^{\omega{>}}\zeta\mbox{ or }\eta
\vartriangleleft
\nu \mbox{ where }\nu\in{\rm supp}(\bools^\alpha_n)
\mbox{ for some }\alpha\in w_{n(*)+1}, n<\omega\},\]
$w'=:\{\alpha<\alpha^*:(\forall n<\omega)(\eta_\alpha\restriction n
\in I)\}$
and $C'=:\{z^*\}$ here; clearly the demands there hold, recalling
${\rm supp}(\bools^\alpha_n)$ is a finite subset of
$\{\rho \in {}^{\omega>}\blackz(\alpha):
\eta_\alpha\rest n \triangleleft \rho\}$ by $\circledast_\alpha(2)$.
So there is a projection $f$ from $\langle
\bool(I', w') \cup \{z^*\}\rangle_{\bool^c_0}$ onto
$\bool(I', w')$, and so by \ref{3.8}(2) clearly $f(z^*)$ is
based on $I'$. As clearly $w'\subseteq\{\alpha:\blackz(\alpha)<\xi\}\cup
w_{n(*)+1}$, we get
\[f(z^*)\in \bool(I', w') \subseteq\langle
\bigcup\limits_{\varepsilon<\xi}
\bool_{[\varepsilon]} \cup \{a_\alpha :\alpha\in w_{n(*)+1}\}
\rangle_{\bool^c_0},\]
So $f(z^\ast)$  belongs to $\bool(I', w')$, which is $\subseteq\langle
\bigcup\limits_{\varepsilon<\xi}\bool_{[\varepsilon]}\cup\{a_\alpha:\alpha\in
w_{n(*)}\}\rangle_{\bool^c_0}$. Also $f(x_{n(*)})=x_{n(*)}$ and
$f(y_{n(*)})=y_{n(*)}$ as $x_{n(*)}$, $y_{n(*)}\in \bool(I', w')$ , so
as $z^*\cap x_{n(*)}=y_{n(*)}$ by the choice of $z^*$, necessarily $f(z^*)
\cap x_{n(*)}=y_{n(*)}$, so by the previous sentence we get a contradiction
to clause (vi) for $n(*)$. So $(**)$ holds.

\par \noindent
\underline{Part V.} We continue the first paragraph of Part IV, and let 
$(**)$ of Part IV hold for $w^*$ and $x^*$.

Let $d_0,\ldots,d_m \in \bool_{[\xi]}$ be such that
\begin{description}
\item[$\boxtimes$]
\begin{enumerate}
\item[(a)] $\bigcup\limits^m_{\ell=0} d_\ell=1 \mbox{ and }$
\item[(b)] $(\forall\ell\leq m)(\forall\alpha\in w^*)
(d_\ell\leq a_\alpha\vee d_\ell\cap a_\alpha=0)$.
\end{enumerate}
\end{description}
There is an $\ell\leq m$ such that $\{y\cap d_\ell:y\leq x^*$ and $y\in
\bool'\}$ is infinite. It is clear (by Part II) that
$\bool''=\langle \bool',d_\ell\rangle_{\bool^c_0}$ is $\aleph_1$--
compact; also $x^*\cap d_\ell\in \bool''\setminus Fi(\bool'')$.

Now, assume that $y\in \bool''$, $y\leq x^*\cap d_\ell$.
Clearly for some $y'\in \bool'$ we have  $y=y'\cap d_\ell$ and 
\wolog\ 
~$y'\leq x^*$. By $(**)$, that is the choice of $w^\ast, x^\ast$ for some $z\in
\langle\bigcup\limits_{\zeta<\xi}\bool_{[\zeta]}
\cup\{a_\alpha:\alpha\in w^*\}
\rangle_{\bool^c_0}$ we have $z\cap x^*=y'$. 
Hence $z\cap (x^*\cap d_\ell)=y$,
and by the choice of $d_\ell$ that is $\boxtimes (b)$ and the choice of $z$,
for some $z'\in\bigcup\limits_{\zeta<\xi}
\bool_{[\zeta]}$, the equation $z'\cap(x^*\cap d_\ell)=z\cap
(x^*\cap d_\ell)=y$
holds.

So by the previous paragraph, in $\bool''$ 
the element $x^{**}\stackrel{\rm
def}{=} x^* \cap d_\ell$ satisfies the requirements in 
$(**)$ for $w^{**} =:
\emptyset$ . Now we use (c) of part I. As $\cf(\xi)=\aleph_0$, let
$\xi=\bigcup\limits_{n<\omega}\zeta_n$ with $\zeta_n<\zeta_{n+1}<\omega$,
and by induction on $n<\omega$ we choose $x_n,y_n$ such that:
\begin{description}
\item[(i)]   $x_n\in \bool''$,  $x_n\leq x^{**}$, and $m<n\ \Rightarrow\ x_m\cap
x_n=0$,
\item[(ii)]  $x^{**}-\bigcup\limits_{\ell<n} x_i\notin Fi(\bool'')$,
\item[(iii)] $y_n\in \bool^{''}$, $y_n\leq x_n$,
\item[(iv)]  for no $z \in \bool_{[\zeta_n]}$ do we have $z\cap x_n=y_n$.
\end{description}
As $\bool''$ is $\aleph_1$--compact, for some $z^*\in \bool''$
we have $z^*\cap x_n= y_n$ for every $n$.

Now, as $\bool''$, $x^{**}$, $w^{**}=\emptyset$ satisfy $(**)$, for some
$z^{**}\in\bigcup\limits_{\zeta<\xi} \bool_{[\zeta]}$ we have $z^*\cap
x^{**}=z^{**}\cap x^{**}$. So for some $n$, $z^{**}\in \bool_{[\zeta_n]}$,
contradicting clause (iv) above. Thus we have finished the proof of
\ref{3.10}. \qquad \hfill$\qed_{\ref{3.10}}$

\begin{Claim}
\label{3.11}
$\bool_{\alpha^*}$ is endo-rigid.
\end{Claim}

Before proving \ref{3.11} we prove the subclaim \ref{3.11A} (For endomorphism
$h$ of $\bool_{\alpha^\ast}$ we shall try to find $\alpha\in Y'$ such that
$h(a_\alpha)$ has to realize $p_\alpha$ to get contradiction, but before choosing $\alpha$
we try to choose appropriate $\langle d^\alpha_n:n<\omega\rangle$, this
is what \ref{3.11A} does for us):

\begin{Subclaim}
\label{3.11A}
Assume that $h$ is an endomorphism of 
$\bool_{\alpha^*}$ and $\bool_{\alpha^*}/
\ExKer(h)$ is an infinite Boolean algebra. \Then\ we can
find $\rho^*$ and $\bar{d}$ such that
\begin{description}
\item[(A)] $\bar{d}=\langle d_n:n<\omega\rangle$ and
$\rho^*\in{}^{\omega{>}}\lambda$,
\item[(B)] $\{d_n:n<\omega\}$ is a maximal antichain of
 $\bool_{\alpha^*}$ and $d_n>0$ of course,
\item[(C)] at least one of $\boxtimes_1$, $\boxtimes_2$,
$\boxtimes_3$ 
hold, 
where
\begin{description}
\item[$\boxtimes_1$]
\begin{description}
\item[(a)] if $\rho^*\vartriangleleft\rho^{**}\in
{}^{\omega{>}}\lambda$  and $n\in (0,\omega)$ then for some  $s\in
\langle x_\nu:\rho^{**}\vartriangleleft\nu\in {}^{\omega{>}}\lambda
\rangle_{\bool^c_0}\setminus\{0, 1\}$ we have $h(s)\cap d_n =0$,
\item[(b)] for no $x\in \bool^c_0$  do we have
$n<\omega\ \Rightarrow\ x\cap h(d_{2n}
\cup d_{2n+1})= h(d_{2n+1})$.
\end{description}
\item[ $\boxtimes_2$] for no $x\in \bool^c_0$ do we have: for every $n\in
(2, \omega)$
\[\begin{array}{rcl}
n\mbox{ is odd}&\Rightarrow&x\cap h(d_n)\cap d_0=h(d_n)\cap d_0,\
\mbox{ and}\\
n\mbox{ is even}&\Rightarrow& x\cap h(d_{n})\cap d_0=0.
  \end{array}\]
\end{description}
\end{description}
\end{Subclaim}

\par \noindent
Proof. As in the proof of \ref{2.4}(2), we can ignore the
maximality requirement in clause (B) (call it (B)$^-$).

Recall
\[\ExKer^*(h)=\{a\in \bool_{\alpha^*}: \{x/ \ExKer(h): x\leq a\}\mbox{ is
finite}\}.\]  
Let  $\idealI_h=\{a\in \bool_{\alpha^*}:$ the set 
$\{h(d):d\leq a$ in $\bool_{\alpha^*}\}$ is finite$\}$, clearly it 
is an ideal of $\bool_{\alpha^*}$ included in $\ExKer^*(h)$ hence $1_{\bool_{\alpha^*}}
\notin \idealI_h$

\noindent {\underline{Case $\alpha$}:}\quad 
For some $\rho^*\in{}^{\omega{>}}
\lambda $ and $a^*\in \bool_{\alpha^*}\setminus \idealI^*(h)$ 
we have: for every
$\rho$ satisfying 
$\rho^*\vartriangleleft\rho\in{}^{\omega{>}}\lambda$ there
is $s\in\langle\{ x_\eta:\rho\vartriangleleft\eta\in
{}^{\omega{>}}\lambda\}
\rangle_{\bool^c_0} \setminus\{0_{\bool^c_0}\}$ such that 
 $h(s\cap a^*)=0$.

\noindent Without loss of generality 
${\rm supp} (a^*)\subseteq \{\rho\in {}^{\omega>}\lambda:
\neg (\rho^*\triangleleft \rho)\}$, hence above $s\cap a^*\neq 0_{\bool^c_0}$ .
Let $\bool^a=\{h(d):d\leq a^*\}\cup\{1-h(d):d\leq a^*\}$. This is a
Boolean subalgebra of $\bool^c_0$ and $1-h(a^*)$
is an atom in it (or zero). As $a^*\notin 
\idealI_h$, clearly $\bool^a$ is infinite, hence by \ref{3.10} there is an
antichain $\langle y_n:n<\omega\rangle$ of $\bool^a$ such that for no $x\in
\bool^c_0$ do we have $x\cap (y_{2n}\cup y_{2n+1})=
y_{2n}$. \Wlog $y_n\leq
h(a^*)$ as at most one $y_n$ fails this. 
Let $d_n$ be such that $h(d_n)=
y_n$ and \wolog\ $d_n\leq a^*$; of course 
$y_n>0$, hence $d_n >0$. 
\Wlog\ $\{d_n:n<\omega\}$ is an antichain (as we can use $d_n-
\bigcup\limits_{m<n} d_m$).

Let $d'_0=1-a^*$, $d'_1=d_0 \cup d_1$ and $d'_{2+n}=d_{2+n}$. 
So clearly clauses (A),
(B)$^-$, (C)$\boxtimes_1$ hold for $\langle d'_n:n<\omega\rangle$, and our
$\rho^*$.
\smallskip

\noindent{\underline{Case $\beta$}:}\quad For some $a^*\in \bool_\alpha$,
$\{h(x)-a^*:x\in \bool_{\alpha^\ast}, x\leq a^*\}$ is infinite.

\noindent Clearly
\[\bool^a=\{h(x)-a^*:x\in \bool_{\alpha^*}\ \&\ x\leq a^*\}
\cup\{1-(h(x)-a^*):x\in
\bool_{\alpha^*}\ \&\ x\leq a^*\}\]
is a subalgebra of $\bool_{\alpha^*}$ (and $a^*$ is an atom in it). By the
assumption (of this case) $\bool^a$ is infinite. So by \ref{3.10} there are
pairwise disjoint 
$y_n\in \bool_a\setminus \{0\}$ such that $\neg(\exists x\in
\bool_a)\bigwedge\limits_{n<\omega}(x\geq y_{2n+1}\ \&\  x\cap 
y_{2n}=0)$. As
$a^*$ is an atom of $\bool^a$, \wolog\ $y_n \leq 1-a^*$, 
hence there are $d_n \in
\bool_{\alpha^*}$ such that $d_n \leq a^*$ and $h(d_n) - a^* =
y_n$. Clearly
\[h(d_n - \bigcup\limits_{\ell< n} d_\ell)= y_n -
\bigcup\limits_{\ell< n} y_\ell = y_n\]
hence \wolog\ the $d_n$ are pairwise disjoint. Let
$d'_0=1-a^*, d'_1= d_0\cup d_1$ and 
$d'_{2+n} = d_{2+n}$, so $\langle d'_n: n<\omega\rangle$ is an
antichain and  
$h(d'_{n}) \cap d_0 = h(d_n) - a^* = \ebool_n$ for $n=2,3,\ldots$, 
hence for no $x\in
\bool_{\alpha^*}$ do we have $n<\omega\Rightarrow 
x\cap h(d'_{2+2n} \cup d'_{2+2n+1}) =
h(d'_{2+2n+1})$. So $\langle d'_n: n<\omega\rangle$ are as requested in
$\boxtimes_2$.
\begin{center}
{\bf Why the two sub-cases exhaust all possibilities?}
\end{center}
Suppose none of Cases $(\alpha)$, $(\beta)$ occurs.
As case ($\alpha$) fail for $a^*=1_{\bool^{\alpha^*}}$ necessarily for 
some $\rho_0\in {}^{\omega>}\lambda$ we have
\begin{enumerate}
\item[(a)] $h(s)>0$ for every $s\in \langle\{x_\eta:\rho_0
\trianglelefteq \eta\in {}^{\omega>}\lambda\}\rangle_{\bool^c_0}
\setminus \{0,1\}$.

Clearly $a\in \langle\{x_\eta:\rho_0\triangleleft \eta\in 
{}^{\omega>}\lambda\{\rangle_{\bool^c_0}\setminus\{0\}$ implies that 
$a\notin \idealI_h$.
As clause $(\beta)$  fail clearly for every $a\in \bool_{\alpha^*}$ the set 
$\{h(d)-a:d\leq a, d\in \bool_{\alpha^*}\}$ is finite. Next we note that:

\par \noindent
$\boxplus$ \ if $\rho_0\triangleleft \rho \in {}^{\omega>}\lambda$ then for some 
$s\in \langle\{x_\eta:\rho \triangleleft \eta \in {}^{\omega>}\lambda\}
\rangle_{\bool_0}\setminus\{0,1\}$ we have $h(s)=s$.

${}$

[Why? for each $\alpha<\omega_1$ let $n_\alpha=(\{h(d)-x_{\rho\conc\langle\alpha\rangle}
:d\leq x_{\rho\conc\langle\alpha\rangle}\})+ |\{h(d)-(-x_{\rho\conc\langle\alpha\rangle}):
d\leq(-x_{\rho\conc\langle\alpha\rangle})\}|$, so we know that $n_\alpha$ is enough, so 
for some $n(*)$ the set $Z=\{\alpha<\omega_1:n_\alpha=n(*)\}$ is infinite. By Ramsey 
theorem for some $a$ we have: if $\alpha<\beta$ are from $Z$ and ${\bf t}_1, {\bf t}_2
\in \{0,1\}$ are truth values then $h(x^{{\bf t}_1}_{\rho\conc\langle\alpha\rangle}\cap
 x^{{\bf t}_2}_{\rho\conc\langle\beta\rangle})-x^{{\bf t}_1}_{\rho\conc\langle\alpha\rangle}
=a_{{\bf t}_1,{\bf t}_2}\in \bool_\alpha$, $h(x^{{\bf t}_1}_{\rho\conc\langle\alpha\rangle}
\cap x^{{\bf t}_2}_{\rho\conc\langle\beta\rangle})\cap x^{{\bf t}_2}_{\rho\conc\langle\beta\rangle}
=b_{{\bf t}_1,{\bf t}_2}\in \bool_\alpha$ where $x^{\bf t}$ is $x$ if ${\bf t}$ is 
1 and is $-x$ if ${\bf t}$ is 0. Let $\alpha_0<\alpha_1<\alpha_2<\alpha_3$ be from $Z$ 
and let $s=x_{\rho\conc\langle\alpha_0\rangle}\cap (-x_{\rho\conc\langle\alpha_1\rangle})\cap
x_{\rho\conc\langle\alpha_2\rangle}\cap (-x_{\rho\conc\langle\alpha_3\rangle})$. 
Now $h(s)\leq x_{\rho\conc\langle\alpha_0\rangle}$ as $h(x_{\rho\conc\langle\alpha_2\rangle})
-x_{\rho\conc\langle\alpha_0\rangle}=h(x_{\rho\conc\langle\alpha_3\rangle})-
x_{\rho\conc\langle\alpha_0\rangle}$ and the equation above, similarly 
$h(s)\leq x_{\rho\conc\langle\alpha_1\rangle}$ and also $h(s)\leq 
x_{\rho\conc\langle\alpha_2\rangle}$ (using $x_{\rho\conc\langle\alpha_0\rangle},
x_{\rho\conc\langle\alpha_1\rangle})$ and $h(s)\leq (-x_{\rho\conc\langle\alpha_3\rangle})$
together $h(s)\leq s$. Now $s_1=s-h(0)>0$. Easily $s>0$ and $s$ is disjoint to 
$b^*=\cup\{a_{{\bf t}_1,{\bf t}_2}\cup b_{{\bf t},{\bf t}_2}:{\bf t}_1,{\bf t}_2$ 
truth values$\}$. 
If $({\bf t}_0,{\bf t}_1)\neq (1,0)$ and $s_1\cap h(x^{{\bf t}_0}_{\rho\conc
\langle\alpha_0\rangle}\cap x^{{\bf t}_1}_{\rho\conc\langle \alpha_1\rangle})
>0$, as $s\cap b^*=0_{\bool_\alpha}$ we get easy contradiction. Similarly for
$(x_{\rho\conc\langle\alpha_2\rangle},x_{\rho\conc\langle\alpha_3\rangle})$ hence 
$h(s)=s$.

So $\boxplus$ holds. 
Let $a^* \in \bool_{\alpha^*}$ be
such that $h(a^*)\nleq a^*$ (exists by \ref{2.4}), and let $a^{**}=
h(a^*)-a^*>0$. By ``not Case $(\alpha)$'' and $\boxplus$, 
for some $\rho_0\in {}^{\omega{>}} \lambda$,

\item[(b)] $h(s)\cap a^{**}\neq 0$ for every 
$s\in\langle x_\eta:\rho_0
\vartriangleleft\eta\in{}^{\omega{>}}\lambda\rangle_{\bool^c_0}\setminus
\{0,1\}$.

Possibly increasing $\rho_0$ the set $\{\eta: \rho_0 \vartriangleleft
\eta\in {}^{\omega>}\lambda\}$ is disjoint to ${\rm supp}(a^*)\cup
{\rm supp}(h(a^*))$. Let $s_n = x_{\rho^0\conc\langle
n\rangle}-\bigcup\limits_{m<n} x_{\rho^0\conc\langle m \rangle})$ for
$n< \omega$, so the $s_n$'s are pairwise 
disjoint non-zero members of
$\bool_{\alpha^*}$ and by (a) we have $h(s_n)\cap a^{**} \neq 0$. 
But
$h(s_n)\cap h(a^*) - a^* = h(s_n)\cap (h(a^*) - a^*) = 
h(s_n) \cap
a^{**} > 0$. So clearly the assumption of case ($\beta$) holds (for $a^\ast$).
\end{enumerate}
 \qquad \hfill$\qed_{\ref{3.11A}}$
\medskip

\underline{Proof of \ref{3.11}.} Suppose $h$ is a counterexample, i.e., $h$
is an endomorphism of $\bool_{\alpha^*}$ but $\bool_{\alpha^*}/
\ExKer(h)$ is
infinite, and we shall get a contradiction.

Clearly if for some good candidate $\alpha$, $h_\alpha\subseteq h$ and
$\alpha\in Y_1$ (see Stage B) then $h(a_\alpha)$ realizes the type $p_\alpha=
\{x\cap b^\alpha_n=c^\alpha_n: n<\omega\}$, a contradiction
(as by clause $\circledast(4)$ of stage A,
$\bool_{\alpha^*}$ omits $p_\alpha$). So we shall try to find such
$\alpha$ which satisfies the requirements $\circledast^1_\alpha$ of 
stage B (hence implicitly  $\circledast_\alpha$ of stage A)
for belonging to
$Y_1$. Let $\rho^*$, $\langle d_n: n<\omega\rangle$ be as in \ref{3.11A}
$(\rho^\ast$ is needed only if $\boxtimes_1$ of (C) of \ref{3.11A} holds otherwise 
we can let
$\rho^\ast=\langle\rangle$) and let $\xi<\lambda$ be such that
\[\bigcup\{{\rm supp}(d_n)\cup{\rm supp}(h(d_n)):n<\omega\}\subseteq
{}^{\omega{>}}\xi\]
Let  $n(*)={\rm lg}(\rho^*)$.
Let $Z\subseteq {}^{\omega>}\lambda$ , it will be used only in 
\ref{3.12AT}. 
We can find a good candidate $\alpha< \alpha^\ast$ such that
\begin{description}
\item[(a)] $h_\alpha\subseteq h$, and $\blackz(\alpha)>\xi$
\item[(b)] $d_n\in \bool[N^\alpha_0]$ for $n<\omega$ and $\rho^*\in\Rang(
f^\alpha)$, and $\langle d_n: n<\omega\rangle \in N^\alpha_0$
\item[(c)] $N^\alpha_0$ is an elementary submodel of the expansion
\\
$({\cal H}_{<\aleph_1}(\lambda), \in, \bool_{\alpha^*}, h, \{(\eta,
x_\eta): \eta\in {}^{\omega>}\lambda\},Z)$ of ${\cal M}$
\end{description}
so in particular  $N^\alpha_n$ is closed under the functions implicit
in the choice of $\bool^\alpha_0$ and $\rho^*$, $\langle d_n:n<\omega\rangle$,
i.e.,
\begin{description}
\item[(d)] $a\in \bool[N^\alpha_0]\ \Rightarrow\ {\rm supp}(a)\subseteq
N^\alpha_0$,
\item[(e)] $\eta\in N^\alpha_m\cap {}^{\omega{>}}\lambda\ 
\Leftrightarrow\
x_\eta\in \bool[N^\alpha_m]$,
\item[(f)] $({}^{\omega{>}}\lambda)\cap N^\alpha_m$ is closed under initial
segments, and each node has infinitely many immediate successors,
\item[(g)] if $\boxtimes_1$ of clause C of \ref{3.11A} holds and if
$\rho^*\vartriangleleft\rho^{**}\in N^\alpha_n$, and $n$ 
is large enough then there is $s\in N^\alpha_n$ as required in
\ref{3.11A}(C)$\boxtimes_1$ so $h(s)-d_0=0$.
\end{description}
As $\blackw$ is a barrier this is possible (using the game
$\Gm'(\blackw)$ and not
$\Gm(\blackw)$ because of the requirement 
$\rho^*\in\Rang(f^\alpha)$ 
recalling Definition \ref{1.7}, that is we 
choose a strategy for player I, 
choosing the $N_n$-s and in the zeroth move also 
$f_\ell$ for $\ell=1,\ldots, \lg (\rho^*)+1$. 
So for some play of the game, player II wins while 
player I uses the strategy described above so the play is 
$\langle (N^\alpha_n,f^\alpha_n): n<\omega\rangle$ for some 
$\alpha<\alpha^*$, so we are done).
Note that the proof of \ref{3.12AT} below use the rest of the present proof, only
ignoring case III below. We then
will choose $\eta_\alpha$, an $\omega$--branch of $\Rang(f^\alpha)$ above
$\rho^*$; but $\blackw$ is a disjoint barrier (see Definition 
\ref{1.6A}(3)) hence $\eta_\alpha$  
hence distinct from $\eta_\beta$ for $\beta<\alpha$
and we will choose 
$s_n\in N^\alpha$ in $\langle x_\nu:\eta_\alpha\restriction
n\vartriangleleft\nu\in {}^{\omega{>}}
\lambda\rangle_{\bool^c_0}\setminus \{0,1\}$
 and let
$b_n=h(d_n)$, $c_n=h(d_n\cap s_n)$ for $n<\omega$, 
$p_\alpha=\{x\cap b_n=c_n:n<\omega\}$,
and $a_\alpha=\bigcup\limits_{ n<\omega}(d_n\cap s_n)\in 
\bool^c_0$. All should
have superscript $\overline{s}$ (where $\overline{s}=
\langle s_n:n<\omega\rangle$), but we usually omit it or write 
$a_\alpha[\overline s], p_\alpha[\overline s]$ etc. 
It is enough to prove that for at least
one such $\bar{s}$ we have 
$a_\alpha [\bar{s}]$, $\bar{s}[d_n:n<\omega]$ exemplify
that $\alpha\in Y_1$.

The choice of $\overline s$ 
(and $\eta_\alpha$ which is determined by
$\overline s$) is done by listing the demands on 
them (see Stage B) and
showing that a solution exists. The only problematic one is (4) (omitting
$p_\beta$ for $\beta\leq\alpha$, $\beta\in Y_1$) and we partition it to
three cases:
\begin{description}
\item[(I)]   $\blackz(\beta)<\blackz(\alpha)$ or $\blackz(\beta)=
\blackz(\alpha)$, $\beta+2^{\aleph_0}\leq \alpha$,
\item[(II)]  $\blackz(\beta)=\blackz(\alpha)$, $\beta<\alpha<\beta+2^{\aleph_0}$,
\item[(III)] $\beta=\alpha$.
\end{description}

We shall prove first that every $\overline s$ is 
O.K. for (I), second that
for any family $\{(\eta^i,\overline s^i): 
i<2^{\aleph_0}\}$ ($\eta^i$ is a
branch of $\Rang(f^\alpha)$ above $\rho^*$, etc.) with pairwise distinct
$\eta^i$'s, all except $<2^{\aleph_0}$ many are O.K. for any instance of
(II), and third that for every $\eta$ 
(a branch of $\Rang(f^\alpha)$ above
$\rho^*$) there is $\overline s$ such that 
$(\overline\eta,\overline s)$
satisfies (III). This clearly suffices (as for each branch $\eta$ 
of Rang$(f^\alpha)$ choose $\overline{s}_f$ \st\ 
$(\eta,\overline{s}_\eta)$ satisfies III, and then chose $\eta$ 
 \st\ $(\eta,\overline{s}_\eta)$ satisfies II) .
\medskip

\noindent\underline{Case I}:\quad $\blackz(\beta)<\blackz(\alpha)$ or $\blackz(
\beta)=\zeta(\alpha)$, $\beta+2^{\aleph_0}\leq\alpha$

\noindent Let $\overline{s}$ be as above. 
Suppose some $x\in\langle \bool_\alpha,a_\alpha[\overline s]
\rangle_{\bool^c_0}$ realize $p_\beta$. Clearly there is a 
partition $\langle
y_\ell:\ell < 4\rangle$ of 1 (in $\bool_\alpha$) 
such that $x=y_0\cup(y_1\cap
a_\alpha [\overline s])\cup (y_2-a_\alpha 
[\overline s])$.  Choose $\xi<
\blackz(\alpha)$ large enough and finite $k<\omega$ so that 
\begin{enumerate}
\item[$\boxdot$] $[\blackz(\beta)<
\blackz(\alpha)\Rightarrow\blackz(\beta)<\xi]$, and $d_n,h_\alpha(d_n),
b^\beta_n$, are based on 
$\{x_\nu:\nu\in {}^{\omega>}\xi\}$ (for $n<\omega$)
and $c^\beta_n$ (for $n<\omega$), $y_0$, $y_1$, 
$y_2$, $y_3$ are based on
$J=\{x_\nu:\nu\in {}^{\omega{>}}\lambda,\ \eta_\alpha\restriction k\not
\vartriangleleft\nu\}$, where $k<\omega$ also satisfies that
$\eta_\alpha(k)>\xi$, $\eta_\alpha\restriction k\notin N^\beta$ (where
$\eta_\alpha\in {}^{\omega}\lambda$ is the one determined by $\overline s$).
\end{enumerate}
These are possible
because of \ref{3.2} and \ref{1.7A}(2)(e).

We claim:
\begin{description}
\item[(*)] there is $m<\omega$ such that 
$b^*=(b^\beta_m\cap (y_1\cup y_2))-
\bigcup\limits_{n\leq k} d_n\neq 0$.
\end{description}
For suppose (*) fails, then as 
$a_\alpha[\overline \bools]\cap (\bigcup\limits_{
n\leq k} d_n)\in \bool_\alpha$; \wolog  
\[(y_1\cup y_2)\cap\bigcup\limits_{n\leq k} d_n=0\]
[Why? otherwise let
\[\begin{array}{l}
y'_0=y_0\cup(y_1\cap a_\alpha
[\overline s]\cap\bigcup\limits_{n\leq k} d_n)
\cup (y_2\cap(\bigcup\limits_{n\leq k} d_n-a_\alpha
[\overline s]))\\
y'_1=y_1-\bigcup\limits_{n\leq k} d_n,\\
y'_2=y_2-\bigcup\limits_{n\leq k} d_n\quad .
  \end{array}\]
So for every $m<\omega$, $b^\beta_m\cap (y'_1\cup y'_2)=0$.
We should now check that the demands on $k$ in $\boxdot_{\xi,k}$ 
are still satisfied (for $\xi$ there is no change)].

Thus, if $x$ realizes $p_\beta$ then so does $y_0$, but $
y_0\in \bool_\alpha$
contradicting the induction hypothesis. So (*) holds.

Now as $\langle d_n:n<\omega\rangle$ is a maximal 
antichain in $\bool_\alpha$, for some $\ell<\omega$,
\[d_\ell\cap b^*=d_\ell\cap(b^\beta_m\cap((y_1\cup y_2)-
\bigcup\limits_{n\leq k} d_n))\neq 0.\]
Necessarily $\ell> k$. So for some $i \in\{1,2\}$ we have $d_\ell\cap
b^* \cap y_i\neq 0$. As $x$ realizes $p_\beta$, necessarily $x\cap(
d_\ell\cap b^\beta_m\cap y_i)=d_\ell\cap c^\beta_n\cap y_i$,
which is based on $J$. But we know that $x\cap 
(d_\ell\cap b^\beta_m\cap y_i)$ is
\[\begin{array}{ll}
d_\ell\cap b^\beta_m\cap y_1\cap a_\alpha [\overline s]=
d_\ell\cap b^\beta_m \cap y_1\cap s_\ell&\mbox{ (if $i=1$)}\\
\mbox{ or } d_\ell\cap b^\beta_m\cap y_2\cap
(1-a_\alpha([\overline s])=d_\ell\cap
b^\beta_m\cap y_2\cap (1-s_\ell)&\mbox{ (if $i=2$).}
  \end{array}\]
As $d_\ell\cap b^\beta_m\cap y_i\neq 0$ is based on $J$, $\ell>k$,
$\eta_\alpha(k)>\xi$, clearly $s_\ell$ is free over 
$\{x_\nu:\nu\in J\}$
(see Fact \ref{3.3}(1)). As $d_\ell \cap b^\beta_m\cap 
y_i\geq d_\ell
\cap b^\ast \cap y_i>0$ and $s_n\notin 0,1$ necessarily
$x\cap (d_\ell\cap b^\beta_m\cap y_i)$
is not based on $J$, contradiction.
\medskip

\noindent\underline{Case II:}\quad $\beta<\alpha<\beta+2^{\aleph_0}$.

\noindent We shall prove that if $\eta^i$, $\bar{s}^i$ 
are appropriate
(for $i=1,2$) and $\eta^1\neq\eta^2$ then $p_\beta$ cannot be realized in
both $\langle \bool_\alpha,a_\alpha[\bar{s}^i]\rangle_{\bool^c_0}$. (So as $\beta<\alpha<\beta+2^{\aleph_0}$, there are less than $2^{\aleph_0}$
non-appropriate pairs ($\eta^i$, $\bar{s}^i$)).

So toward contradiction, for 
$i=1,2$, let $x^i\in\langle \bool_\alpha,
a_\alpha[\overline{s}^i]
\rangle_{\bool^c_0}$ realize $p_\beta$. 
Clearly there is a partition $\langle
y^i_\ell:\ell<4\rangle$ of 1 (in $\bool_\alpha$) such that
\[x^i=y^i_0\cup (y^i_1\cap a_\alpha[\overline{s}^i])
\cup (y^i_2-a_\alpha[\overline{s}^i]).\]
Choose $\xi<\blackz(\alpha)$ large enough and finite $k<\omega$ such that
\begin{description}
\item[(i)]   $d_n,h_\alpha(d_n), b^\beta_n$ (for $n<\omega$) 
are based on $\{x_\eta:\eta\in
{}^{\omega{>}}\xi\}$,
\item[(ii)]  $y^i_\ell$ (for $i=1,2$ and
$\ell<4$) and $c^\beta_n$ (for $n<\omega$) are based on
\[J=\{x_\nu:\nu\in {}^{\omega{>}}\lambda\ \&\ \eta^1
\rest k \not\vartriangleleft\nu\
\&\ \eta^2\rest k\not\vartriangleleft\nu\},\]
\item[(iii)] $\eta^1(k)>\xi$, $\eta^2(k)>\xi$ and $\eta^1\restriction k\neq
\eta^2\restriction k$.
\end{description}
We claim that
\begin{description}
\item[$(*)$] there is $m<\omega$ such that $0<b^*=:
b^\beta_m-(y_0^1\cup y_3^1
\cup y^2_0\cup y^2_3)-\bigcup\limits_{n\leq k} d_n$.
\end{description}
[Why? Otherwise $a^i=:a_\alpha[\overline{\bools}^i]\cap
(y^i_0\cup y^i_3\cup
\bigcup\limits_{n\leq k} d_n)$ belongs to $\bool_\alpha$ for $i=1,2$ and
$a^1\cup a^2$ realizes $p_\beta$, a contradiction.]

Clearly $b^*$ is based on $J$.

As $\langle d_n:n<\omega\rangle$ is a maximal antichain in $\bool^c_0$ (and
hence in $\bool_\alpha$), for some $\ell<\omega$ we have $0<d_\ell\cap
b^*$; clearly $\ell>k$. So for some $j(1),j(2)\in\{1,
2\}$ we have $0<b^{**}=:d_\ell\cap b^* \cap y^1_{j(1)}\cap
y^2_{j(2)}$ (just recall $y^i_1\cup y^i_2=1-
(y^i_0\cup y^i_3)$).
Clearly also $b^{**}$ is based on $J$ and $b^{**}\leq d_\ell
\cap b^*\leq d_\ell\cap b^\beta_m$ by the choice of $b^{**}, b^*$ 
respectively. So by the last two sentences, 
as $x^i$ realizes $p_\beta$, clearly
$x^i\cap (d_\ell\cap b^\beta_m)=d_\ell\cap c^\beta_m$,
but the latter does not depend on $i$. Hence $x^1\cap
b^{**}=x^2\cap b^{**}$. But as $\overline{a}_\alpha
[\overline{s}^i]=\bigcup\limits_n (d_n \cap s^i_n)$ 
we know
that $x^i\cap d_\ell$ is $d_\ell\cap s^i_2$ if $j(i)=1$ and is
$d_\ell- s^i_\ell$ if $j(i)=2$. We can conclude that
\[\mbox{either }\quad b^{**}\cap s^1_\ell=b^{**}\cap 
s^2_\ell\quad\mbox{ or }\quad b^{**}\cap 
s^1_\ell=b^{**}\cap (-s^2_\ell)\]
(the other 2 possibilities are reduced to those). 
But $b^{**}$ is based on
$J$ whereas ${\rm supp}(s^1_\ell)$, ${\rm supp}
(s^2_\ell)$ are disjoint
subsets of $\{x_\eta:\eta\in {}^{\omega{>}}\lambda\}\setminus J$ and
$0<s^i_\ell<1$, a contradiction.
\medskip

\noindent\underline{Case III:}\quad $\beta=\alpha$.

\noindent This case is splitted into two sub-cases. Let $\eta_\alpha$ be
any $\omega$-branch of $f^\alpha$ \st\ 
$\rho^*\triangleleft \eta_\alpha$, so necessarily 
 $\eta_\alpha\neq \eta_\beta$ whenever $\beta<\alpha$. 
The proof splits to cases according which
$\ell\in\{1,2\}$ is such that $\boxtimes_\ell$ of \ref{3.11A}(2) is
satisfied by $\rho^*$,$\bar{d}$. Choose $\rho^*_n\in N^\alpha\cap
({}^{\omega{>}}\lambda)$ such that
\[\ell g(\rho^*_n)=\ell g(\rho^*)+n+1,\quad \rho^*_n\restriction (\ell
g(\rho^*_n)-1)\vartriangleleft\eta_\alpha\quad\mbox{ and}\quad \rho^*_n
\not\vartriangleleft\eta_\alpha.\]
\smallskip

\noindent \underline{Sub-case 1:}\quad $\boxtimes_1$ holds.

\noindent We choose $s_n\in \bool[N^\alpha]$, satisfying 
$s_n\in \langle x_\nu:\rho^*_n
\vartriangleleft\nu\in {}^{\omega{>}}\lambda\rangle_{\bool^c_0}$, 
$s_n\neq 0,1$
and
\[\boxplus\quad n=2m\ \Rightarrow\ h(s_n)=0,\qquad n=2m+1\ \Rightarrow\ 
h(s_n)=1\]
this means $x\cap h(d_{2n}\cup d_{2n+1})=h(d_{2n})$
(this is possible by $\boxtimes_1$ applied to 
$\rho^*_n$ and using $s_n$ or
$1-s_n$).

Assume toward contradiction that $x\in\langle \bool_\alpha\cup
a_\alpha[\bar{s}]
\rangle_{\bool^c_0}$ satisfies $x\cap h(d_n)=h(d_n\cap s_n)$
 hence $x\cap h (d_{2n}\cup d_{2n+1})= h(d_{2n})$
for $n<\omega$. Let $\langle y_\ell:\ell<4\rangle$
be a partition of $1$ in $\bool_\alpha$ such that
$x=y_0\cup (y_1\cap a_\alpha[\bar{s}])
\cup (y_2-a_\alpha[\bar{s}])$.
As the type $q=\{x'\cap h(d_{2n}\cup d_{2n+1})= h(d_{2n}):n<\omega\}$
is not realized in $\bool_\alpha$, and $\langle y_0,
y_1,y_2,y_3\rangle$
is a partition of 1 in $\bool_\alpha$ clearly for some $i<4$ the type
\[q_i=\{x'\cap h(d_{2n}\cup d_{2n+1}) \cap y_i=h
(d_{2n+1})\cap y_i:n<\omega\}\]
Let $k(*)<\omega$ be \st\ $\{\eta:\eta_\alpha\rest k(*)\trianglelefteq \eta
\in {}^{\omega>}\lambda\}$ is disjoint to 
${\rm supp}(y_i)$ for $i<4$ and $\xi<\eta_\alpha(k(*)-1)$
is not realized in $\bool_\alpha$. By the choice of the $y_\ell$'s
and the choice of $x$, necessarily $i\in\{1,2\}$ and for notational 
simplicity let $i=1$.
So ${\cal U}=:\{n:h(d_{2n}\cup h_{2n+1})\cap y_1\cdot>0\}$ is infinite, and 
as $\langle d_k:k<\omega\rangle$ is a maximal antichain of 
$\bool^c_0$ hence of $\bool_\alpha$, clearly
for each $n\in {\cal U}$, the set 
\[{\cal U}_n=\{k:h(d_{2n}\cup d_{2n+1})
\cap d_k\cap y_{\dot{1}}>0\}\]  is nonempty. 
 Clearly  $k\in U_n\Rightarrow x\cap h (d_{2n}\cup d_{2n+1})\cap 
d_k\cap y_1= s_k\cap h (d_{2n}\cup d_{2n+1})\cap d_k\cap y_i$ 
hence $k(*)\leq k\in U_n\Rightarrow {\rm supp} (s_k)\subseteq 
{\rm supp}((x\cap h(d_{2n})\cup d_{2n+1}) \cap y_i)$. Hence if 
$n\in {\cal U}$ and ${\cal U}_n$ is
infinite then $x\cap h(d_{2n}\cup d_{2n+1})\cap y_1$ is not in $\bool_\alpha$,
easy contradiction as $h(d_{2n})\in \bool_\alpha$; so $n\in {\cal U}\Rightarrow
1\leq |{\cal U}_n|<\aleph_0$. If $\cup\{{\cal U}_m:n<\omega\}$ is
finite then $d^*=\cup \{d_\ell:\ell\in {\cal U}_n$ for some 
$n<\omega$ (so $n\in {\cal U})\}$ belong to $\bool_\alpha$ 
(as a finite union of members) and $n<\omega\Rightarrow h(d_{2n}\cup
d_{2n+1})\leq d^*$ and $x\cap d^*\in \bool_\alpha$ so $q_i$ is realized 
and we get easy contradiction. Let $f:{\cal U}
\rightarrow \omega, f(n)= \max ({\cal U}_n)$.
Recall that $k(*)<\omega$ and  
$\xi<\blackz(\alpha)$ are large enough.
For $n\in {\cal U}$ with $f(n)\geq k(*)$, clearly.
\[x\cap (h(d_{2n}\cup d_{2n+1}) \cap y_i\cap d_{f(n)}\in \{h(d_{2n}
\cup d_{2n+1})\cap y_i\cap s_{f(n)},\]
$$
 h(d_{2n}\cup d_{2n+1}) \cap y_i
\cap(-s_{f(n)})\}
$$
by the choice of the $\bar{a}_\alpha[\bar{s}]$'s and
\[x\cap h (d_{2n}\cup d_{2n+1})\cap y_i=h(d_{2n})\cap y_i\]
by the choice of $x$. But the latter, $h(d_{2n})\cap y_i$ 
is supported by $\{x_\nu:\rho^*_n
\not\vartriangleleft\nu\in {}^{\omega{>}}\lambda\}$ 
(as the $\rho^*_m$'s are
pairwise $\vartriangleleft$--incomparable), whereas the former is not by the
choice of $\xi$, and $k(*)$.
\smallskip

\noindent\underline{Sub-case 2:}\quad $\boxtimes_2$ holds.

\noindent We choose $s_n\in\langle x_\nu:\rho^*_n
\vartriangleleft\nu\in
{}^{\omega{>}}\lambda\rangle_{\bool^c_0}$, 
$s_n\in \bool [N^\alpha]\setminus\{
0_{\bool_\alpha},1_{\bool_\alpha}\}$. 
Now for $i\in\{1,2,3,4\}$ we let $\bar s^i=\langle s^i_n:n<\omega\rangle$
be defined as follows: $s^i_n$ is $s_n$ if $n+i$ is even and $-s^i_n$ if 
$n+i$ is odd. If for some $i$ the Boolean algebra 
$\langle\bool_\alpha\cup \bar a_\alpha[\bar s^i]\rangle_{\bool^c_0}$ omit
the type $p^i_\alpha=\{x\cap h(d_n)\cap d_0=h(s^i_n)\cap h(d_n)\cap d_0:n<\omega\}$
then we are done, so assume that $x^i\in \langle \bool_i\cup\{
a_\alpha[\bar s^i]\}\rangle_{\bool^c_0}$ realizes $p^i_\alpha$, hence 
$y_i=: x^i\cap d_0$ realizes $p_\alpha^i$ and belong to $\bool_\alpha$. But then 
$y_1$ realizes $\{x\cap h (d_{2n}\cup d_{2n+1}) \cap d_0=h(d_{2n+1})\cap d_0:n
<\omega\}$  
but this contradict the choice of
$\langle d_n:n<\omega\rangle$, see $\boxtimes_2$ of \ref{3.11A}, so we are done.

So we finish the proof of \ref{3.11}; so $\bool_{\alpha^*}$ is endo-rigid.
\qquad\hfill$\qed_{\ref{3.11}}$

\begin{claim}
\label{3.12AT}
There are $\lambda^{\aleph_0}$ ordinal  
$\alpha<\alpha^\ast$ which belongs to $Y'$ (even to $Y_2$).
\end{claim}
\medskip

Proof: Let $h$ be the identity on $\bool_{\alpha^*}$. In the proof
 of \ref{3.11}, guessing the good candidate $\alpha$ 
we have $\lambda^{\aleph_0}$ possible choices as 
$Z\subseteq {}^{\omega>}\lambda$ was arbitrary and we could use $Z=\{\eta
\rest n:n<\omega\}$ for any $\eta\in{}^\omega\lambda$. We can find
a maximal antichain $\langle d_n:n<\omega\rangle$ of $\bool_0$ included in
$\langle\{x_{\langle \gamma\rangle}:\gamma<\lambda$ moreover
$<\gamma>\in N^\alpha_0\}\rangle_{\bool_0}$. For any $\eta\in \lim (f^\alpha)$
and $\bar{s}=\langle s_n:n<\omega\rangle$ where 
$s_n\in\langle \{x_\rho:
(\eta\rest n)^\wedge \langle \eta(n)+1\rangle\triangleleft \rho\}
\rangle_{\bool^c_0}$ we can define $\bar{a}_\alpha[\bar{s}]=\cup\{
d_n\cap s_n:n<\omega\}\in \bool^c_0$. 
As in the proof of \ref{3.11} for
some such $(\eta,\bar{s})$ the fitting $\eta_\alpha=\eta$,
$a_\alpha=\bar{a}_\alpha[\bar{s}]$, all the demands 
for $\alpha\in Y_2$ (see $\otimes$ in Stage B) are satisfied.
\medskip

\begin{Lemma}
\label{3.12}
$\bool_{\alpha^*}$ is indecomposable.
\end{Lemma}
\medskip

Proof:  Suppose $\idealJ_0,\idealJ_1$ are disjoint ideals of
$\bool_{\alpha^*}$, each with no maximal member,
which generate a maximal ideal
of $\bool_{\alpha^*}$. For $\ell=0,1$ let $\{d^\ell_n:\ell<\omega\}$
be a maximal antichain $\subseteq \idealJ_\ell\setminus
\{0_{\bool_{\alpha^\ast}}\}$ (maximal as subset of
${\cal J}_\ell\setminus\{0_{\bool_{\alpha^\ast}}\})$, 
they are countable as
$\bool_{\alpha^*}$ satisfies the c.c.c., and may be chosen infinite as
$\ell<2\Rightarrow \idealJ_\ell\neq\{0\}$ (and 
$\bool_{\alpha^*}$ is atomless). Let ${\mathcal
J}$ be the ideal generated by $\idealJ_0\cup \idealJ_1$.

Now, for example for some $\xi<\lambda$, $\{d^\ell_n:\ell<2,\ n<\omega\}
\subseteq \bool_{[\xi]}$, so easily for some $\alpha\in Y_2$,
$\blackz(\alpha)>\xi$. Clearly
$a_{\alpha}\in \idealJ$ or $1-a_{\alpha}\in \idealJ$. For
notational simplicity assume $a_\alpha\in \idealJ$. So $a_\alpha=b^0
\cup b^1$, $b^\ell\in \idealJ_\ell$. 
Now, $\rpr_\xi(b^\ell)\in \bool_{[\xi]}$
and is disjoint to each $d^{1-\ell}_n$ so by the maximality of $\{d^{1-
\ell}_n:\ n<\omega\}$, $\rpr_\xi(b^\ell)$ is disjoint to every member of
$\idealJ_{1-\ell}$. As $\idealJ_0\cup \idealJ_1$ generates
a maximal ideal of $\bool_{\alpha^*}$, 
clearly $\rpr_\xi(b^\ell)\in \idealJ_\ell$ [otherwise
$\rpr_\xi(b^\ell)=1-c^0\cup c^1$, 
for some $c^0\in \idealJ_1$, $c^1\in
\idealJ_1$, and then $c^{1-\ell}$ is necessarily a maximal member of
$\idealJ_{1-\ell}$, so $\idealJ_{1-\ell}$ is principal, a
contradiction]. So $\rpr_\xi(b^0)\cup \rpr_\xi(b^1)<1$, 
but $1=\rpr_\xi(a_\alpha)
=\bigcup\limits^2_{\ell=0} \rpr_\xi(b^\ell)$, a contradiction.
\qquad$\qed_{\ref{3.12}}$\ \ $\qed_{\ref{3.1}}$
\medskip

Of course,

\begin{Claim}
\label{3.13}
If $\cf(\lambda)>\aleph_0$ \then\  there are
$\bool_i$ for $i<2^{\lambda^{\aleph_0}}$ such that
\begin{enumerate}
\item[(a)] $\bool_i$ is Boolean algebra of 
cardinality $\lambda^{\aleph_0}$,
density character $\lambda$, and this holds even for
$\bool_i\restriction a$ for
$a\in \bool_i\setminus\{0_{\bool_i}\}$,
\item[(b)] $\bool_i$ is endo-rigid indecomposable,
\item[(c)] any homomorphism from any $\bool_i$ to
$\bool_j$ ($i\neq j$) has a finite
range.
\end{enumerate}
\end{Claim}

Proof: We can repeat the proof of 
\ref{3.1}. Now we build  $\bool_\alpha\langle Z\rangle$
for every $Z\subseteq {}^\omega\lambda$, \st\ for each
$\alpha$ we try to guess
not $\bool_{\alpha^\ast}$ 
and an endomorphism of it but we try to guess
$\bool^1[N^\alpha]=(\bool_\alpha\langle Z_1\rangle)\rest
N^\alpha$, $\bool^2 [N^\alpha]=(\bool_\alpha\langle Z_2\rangle)\rest 
N^\alpha$ and $h=H^{N^\alpha}$ 
an homomorphism from $\bool^1[N^\alpha]$ into $\bool^2
[N^\alpha]$, and we ``kill'' i.e., guarantee $h$ cannot be extended to a 
homomorphism from $\bool\langle Z^1\rangle$, to $\bool\langle Z^2\rangle$
 when $\bool\langle Z^\ell\rangle\rest N^\alpha=\bool^\ell[N^\alpha]$  
 \qquad\hfill$\qed_{\ref{3.13}}$

\begin{Claim}
\label{3.14}
In \ref{3.1}, \ref{3.13} we can replace the assumption
$\cf(\lambda)>\aleph_0$ by $\lambda>\aleph_0$.
\end{Claim}

 Proof: We replace ${}^{\omega>}\lambda,S$ by ${}^{\omega>}
(\lambda\times \omega_1),\{\lambda\times \delta:\delta<\omega_1$ 
a limits ordinal\} so we use 
\cite[3.17]{Sh:309}  
instead of  
\cite[3.11, 3.16]{Sh:309}. 
$\qed_{\ref{3.14}}$

Also note:

\begin{Observation}
\label{3.15}
Assume $2^{\aleph_0}<\lambda<\lambda^{\aleph_0}$, 
and $\bool$ is c.c.c.~Boolean
algebra of cardinality $\lambda$, and there is $\mu, \mu<\lambda<
\mu^{\aleph_0}$, hence \wolog\ $\lambda>\mu=\min\{\mu:\mu^{\aleph_0}
\geq\lambda\}$.
\begin{enumerate}
\item There is a free Boolean algebra $\bool_x$ of
cardinality $\mu$ such that $\bool_0\subseteq \bool$.
\item There is $\bar{\bool}$ such that
\begin{enumerate}
\item[(a)] $\bar{\bool}=\langle \bool_n:n<\omega\rangle$,
\item[(b)] $\bool_n$ is a Boolean subalgebra of $\bool$,
\item[(c)] $\bool_n\subseteq \bool_{n+1}$ and $\bool
=\bigcup\limits_{n<\omega} \bool_n$,
\item[(d)] there is $A_n\subseteq \bool_{n+1}$ independent over
$\bool_n${}\footnote{i.e., for every $a\in \bool_n\setminus
\{0_{\bool_n}\}$ and a non-trivial Boolean combination $b$
of members of $A_n$ we have $a\cap b>0$}
of cardinality $\mu$.
\end{enumerate}
\item $\bool$ is not endo-rigid.
\item There are projections\footnote{i.e homomorphism $h$ from $\bool$ 
into $\bool$ such that $x\in \bool\Rightarrow h(h(x))=h(x)$}
 of $\bool$ whose range are atomless countable
Boolean algebra.
\item there are $\lambda^{\aleph_0}$ atomless Boolean subalgebras
$\bool'$ of $\bool$ \st\ there is a projection from $\bool$ onto $\bool'$
\end{enumerate}
\end{Observation}

Proof: By cardinal arithmetic $(\forall\kappa<\mu)(\kappa^{\aleph_0}<\mu)$ and $\cf(\mu)=
\aleph_0$. Let $\mu=\sum\limits_{n<\omega}\mu_n$, $\mu_n<\mu_{n+1}$,
$\mu_n^{\aleph_0}=\mu_n$.
\medskip

\noindent (1), (2)\quad By \cite[Lemma 4.9, p.88]{Sh:92}, we can find $\langle
b_\alpha:\alpha<\mu\rangle$ independent in $\bool$.
Let $\bool_*$ be the subalgebra
of $\bool$ generated by $\{b_\alpha:\alpha<\mu\}$, and let $\bool^c_*$ be the
completion of $\bool_*$. Let $h^*$ be a homomorphism from $\bool$
into $\bool^c_*$ extending $\id_{\bool_*}$, (it is well known that such 
homomorphism exists) and let $\bool'=\Rang(h^*)$,
so $\bool_*\subseteq \bool'\subseteq
\bool^c_*$, $|\bool_*|\leq |\bool'|\leq\lambda$.
For each $a\in \bool'$ there is a countable
$u_a\subseteq\mu$ such that $a$ is based on 
(i.e., belongs to the completion inside $\bool^c_*$
of) the set $\{b_\alpha:\alpha\in u_a\}$.

We can find pairwise distinct 
$\eta_\alpha\in\prod\limits_{n<\omega}\mu^+_n$
for $\alpha<\lambda^{\aleph_0}$ such that $\eta_\alpha\rest (n+1)\neq
\eta_\beta\rest (n+1)\Rightarrow
\eta_\alpha(n)+\mu_n\neq \eta_\beta(n)+\mu_n$.
Now for each $a\in \bool'$ the set
\[w_a=\{\alpha:(\exists^\infty n)(u_a\cap 
[\mu_n\times\eta_\alpha(n),\mu_n
\times\eta_\alpha(n)+\mu_n)\neq\emptyset)\}\]
has cardinality $\leq 2^{\aleph_0}$. But $|\bool'|+2^{\aleph_0}\leq\lambda<
\lambda^{\aleph_0}$, hence for some $\alpha<\lambda^{\aleph_0}$ we have
\[\eta_\alpha\notin \bigcup\{w_a:a\in \bool'\}\]
Let
\[\begin{array}{ll}
\bool_m'=:&\big\{x\in \bool':h^*(x)\mbox{ is based on } A'_n\}
\end{array}\]
 where 
$$
A'_m=:\{b_\beta:\mbox{ if }
n\geq m \mbox{ then }\beta\notin [\mu_n\times
\eta_\alpha(n),\mu_n\times\eta_\alpha(n)+\mu_n)\}\big\}.
$$
The sequence $\langle \bool_n':n<\omega\rangle$ 
is as required except that in
clause (d) if we naturally let $A''_n=\{b_\beta:\beta\in 
[\mu_n\times \eta_\alpha(n), \mu_n\times \eta_\alpha(n)+\mu_n)\}$ 
we get $|A''_n|\geq \mu_n$ instead $|A''_n|\geq \mu$.
So let  $ \omega $ 
be the disjoint union of the
infinite sets $v_n$ for $n<\omega$, and let $\bool_m
=\{x\in \bool: h^*(x)\mbox{ is based on } A_m\}$, where
\[\begin{array}{ll}
A_m=\{ 
b_\beta:& \mbox{ if }n<\omega\ \mbox{ and }\\
&n\notin\bigcup_{k\leq m}v_k\ \mbox{ then } \beta\notin [\mu_n\times
\eta_\alpha(n),\mu_n\times\eta_\alpha(n)+\mu_n)\}.
\end{array}\]
Then the sequence $\langle \bool_m:m<\omega\rangle$ is as required.
\medskip

\noindent (3) Follows by (4). 

\noindent (4) Choose $a_n\in A_n$ for $n<\omega$. 
Now we define by induction on $n$, a 
projection  $ h_n $ 
 from the Boolean 
algebra $\bool_n$ onto the subalgebra $\bool^*_n$ of $\bool_n$ 
generated by $\{a_\ell:\ell<n\}$ 
freely 
and extending $h_m$ for $m<n$.
For $n=0$, let $D_0$ be any ultrafilter of $\bool_0$ and let 
$h_0 (x)$ be $1_{\bool_0}=1_{\bool}$ if $x\in D_0$ and $0_{\bool_0}
=0_\bool$ if $x\in \bool_0\setminus D_0$. For $n=m+1$ let 
$\langle a^m_k:k<2^m\rangle$ list the atoms of $\bool^*_m$, which is a 
finite Boolean algebra, and for $k<2^m$ let $D^m_k=\{x\in \bool_m:
a^m_k\subseteq h(x) \in \bool^*_m\}$, this is an ultrafilter of
$\bool_m$. For each $k$ we can find two ultrafilters $D^m_{k,0},
D^m_{k,1}$ of $\bool_n=\bool_{m+1}$ extending $D^m_k$ \st\ 
$a_m\in D^m_{k,1}$ and $a_m\notin D^m_{k,0}$. Lastly define 
$h_n=h_{m+1}:\bool_n\rightarrow \bool^*_n$ by $h_n(x)=\bigcup \{
a^m_k\cap a_m:x\in D^m_{k,1}\}\cup \bigcup \{a^m_k - 
 a_m:
x\in D^m_{k,0}\}$, it is easy to check that $h_n$ is a 
homomorphism from $\bool_n$ onto $\bool^*_n$ and is the identity on 
$\bool^*_n$ and extend $h_m$. 

Clearly $h=\cup\{h_n:n<\omega\}$ is a projection of $\bool=\cup\{
\bool_n:n<\omega\}$ onto $\bool^*=\cup \{\bool^*_n:n<\omega\}$, so 
$h, \bool^*$ are as required.

\noindent (5) By the proof of part (4), that is the arbitrary choice of 
$\langle a_n:n<\omega\rangle \in \prod\limits_{n<\omega} A_n$.
\qquad\hfill$\qed_{\ref{3.15}}$
\medskip
  
\begin{Discussion}
\label{3.16AT}
\begin{enumerate}
\item In \ref{3.15} the only use of the c.c.c. is to find a free 
subalgebra of $\bool$ of cardinality $\mu$. 
\item What about $|\bool|<2^{\aleph_0}$? S.Koppelberg 
and the author noted (independently)
that under ${\rm MA}$  (or just ${\mathfrak p}=2^{\aleph_0}$)
such Boolean algebras are not endo-rigid. Why? let $a_n\in\bool\setminus
\{0_\bool\}$ be pairwise disjoint, let $D_n$ be an ultrafilter of $\bool$ to
which $a_n$ belong, and for $x\in \bool$ let ${\cal U}_x=\{n:x\in D_n\}$.
By ${\rm MA}$ there is an infinite ${\cal U} \subseteq \omega$ \st\ for every
$x\in {\cal U}$ the set ${\cal U}\cap {\cal U}_x$ is finite or the set
${\cal U}\setminus{\cal U}_x$ is finite. Let $h\in \Ext(\bool)$ be $h(x)
=\cup\{a_n:n\in {\cal U}_x\}$ if
${\cal U}\cap {\cal U}_x$ is finite and
$h(x)=1_\bool-\cup \{a_n:n\in{\cal U}\setminus {\cal U}_x\}$ if ${\cal U}
\setminus{\cal U}_x$ is finite.
\item Assume $\mu=\sum\{\mu_n:n<\omega\}$, 
$\mu^\kappa_n=\mu_n<\mu_{n+1}$.
If $\bool$ is a Boolean algebra satisfying the $\kappa^+$- c.c. \st\
$\mu<|\bool|<\mu^{\aleph_0}$ then the 
construction of \ref{3.15} holds. The proof is similar.
\end{enumerate}
\end{Discussion}

\begin{Discussion}
\label{4.20}
We may wonder whether Claims \ref{3.9}, \ref{3.10} can be improved to: if
$d_n\in \bool_{\alpha^*}$ (for $n<\omega$) are pairwise disjoint
non--zero, then for some
$w\subseteq\omega$ there is no $x\in \bool_{\alpha^*}$ satisfying
\[[n\in w\ \ \Rightarrow\ \ x\cap d_n=d_n]\quad\mbox{ and }\quad [n\in\omega
\setminus w\ \ \Rightarrow\ \ x\cap d_n=0].\]
The problem is that $\{\eta_\alpha:\alpha<\alpha^*\}\subseteq {}^\omega
\lambda$ may contain a perfect set and if we are not careful about the
$s^\alpha_n$-s mentioned above we may fail. If $\lambda=\mu^+$,
$\mu^{\aleph_0}=\mu$, then we may try to demand, for each $\zeta^*<\lambda$,
that
\[\langle \bigcup\limits_{n<\omega}{\rm supp}(s^\alpha_n):\alpha<
\alpha^*,\zeta(\alpha)=\zeta^*\rangle\]
is a sequence of pairwise disjoint sets. Alternatively we may look for a
thinner black box (of course, preferably more then just no perfect set of
$\eta_\alpha$'s), 
 see \cite[\S3]{Sh:309}.  
\end{Discussion}

\end{document}